\documentclass[a4paper,10pt]{article}
\usepackage{graphicx}

\usepackage{geometry}

\usepackage[ngerman]{babel}
\usepackage{amssymb,amsmath}
\usepackage[usenames,dvipsnames]{color}
\usepackage{amsfonts}
\usepackage{tikz}
\usepackage{hyperref}
\hypersetup{colorlinks,citecolor=black,filecolor=black,linkcolor=black,urlcolor=black}

\usepackage{framed}

\usepackage{latexsym}
\usepackage{amsthm}
\usepackage{bm}

\usepackage{amscd}
\usepackage[all]{xy}

\renewcommand{\abstract}{\textbf{Abstract}}
\usepackage{cite}
\newcommand{\I}{\mathrm{i}}
\renewcommand{\d}{{\rm d}}
\newcommand{\Di}{\mathbf{d}}
\newcommand{\Mu}{\mathbf{m}}
\newcommand{\Co}{\mathbf{c}}
\newcommand{\Def}{\hbar}
\newcommand{\Ps}{Z}

\DeclareMathOperator{\id}{id}
\newcommand{\Conv}{\mathop{\scalebox{1.4}{\raisebox{-0.2ex}{$\ast$}}}}

\newtheorem{satz}{Satz}[subsection]
\newtheorem{theorem}[satz]{\textit{Theorem}}
\newtheorem{lemma}[satz]{\textit{Lemma}}

\newtheorem{proposition}[satz]{\textit{Proposition}}
\newtheorem{korollar}[satz]{\textit{Corollary}}
\newtheorem{definition}[satz]{\textit{Definition}}
\numberwithin{equation}{section}
\date{}

\vspace{-0.6cm}

\title{\textbf{\textit{On self-adjointness of Poisson summation}}\vspace{-0.4cm}}
\author{\textit{Johannes L\"offler}}
\begin{document}

\maketitle
\vspace{-0.4cm}
\begin{abstract} 
We show that a combination of well-known operators, namely $\I{\tau}\circ{H}\circ\Ps$ is self-adjoint and {\em ad-hoc} related to the $\zeta$ function. Here ${\tau}$ is an involution appearing in Weil's positivity criteria needed for symmetrization, $H$ a regularization operator introduced by Connes \cite{Co2} and $\Ps$ essentially Poisson summation. We elaborate on the Hilbert-P\'olya conjecture, discuss why the Hermite-Biehler theorem, uncertainty relations and cohomologies are interesting in our scenario. \vspace{0.27cm}\end{abstract}
\vspace{-0.4cm}
\subsection*{Introduction}
We hope our exposition contains at least for non-experts new informations and maybe can stimulate other approaches, but we admit that the meaning of some things is still a puzzle.

The zeta function $\zeta$ is  for $\Re(s)>1$ defined by $\zeta(s):=\prod_{p\in\mathbb{P}}\frac{1}{1-p^{-s}}=\sum_{n=1}^{\infty}n^{-s}$ where $\mathbb{P}$ is the set of primes. One way to continue $\zeta$ to the complex plane is the following \cite{Ti}:

Let the Mellin transform of $f$ as usual be defined by
${\mathcal{M}[f](s):=\int_0^\infty\frac{\\\d{t}}{t}t^s{f}(t)}$. For example we will in the following considerations sometimes use the so-called $\Gamma$ function for $\Re(s)>0$ defined as the Mellin transform of $e^{-t}$, for instance we set
${\Gamma(s):=\mathcal{M}[\exp(-t)](s)}$. The $\Gamma$ function admits its continuation by the faculty equation
$\Gamma(1+s)=s\Gamma(s)$.

Let the Fourier transform sending $f$ to $\widehat{f}$ as usual be defined by
$\widehat{f}(p):=\int_{-\infty}^\infty\;\d{x}\;f(x)e^{-2\pi\I{x}p}$.
Mellin transform and Fourier transform are related by the equation
$\mathcal{M}[f](\I{s})=\widehat{({f\circ{e}^{-\id}})}(s)$.

For $\beta\in\mathbb{R}^+\setminus0$ we define the dilation operator $\Di_\beta$ by
${(\Di_\beta{f})(y):=f(\beta{y})}$ also called scaling operator. We have $[\Di_{\beta_1},\Di_{\beta_2}]=0$ and dilations also commute with the Zeta operators
\begin{equation}\Ps^\alpha:=\sum_{n=1}^\infty\Di_{n^\alpha}\end{equation}

With this notation we define $\Psi:\mathbb{R}^{+}\rightarrow\mathbb{R}^{+}$ as usual by
${\Psi(t):=\Ps^2\bigr(\exp(-\pi{t})\bigr)=\sum_{n=1}^\infty{e}^{-\pi{n}^2t}}$. Notice that we can estimate $\Psi$ by a geometric series for instance we have the inequality
$\Psi(t)\leq{\exp(-\pi{t})}/({1-\exp(-\pi{t})})$
and hence $\Psi(t)$ converges exponentially fast to zero as $t\rightarrow+\infty$.

By $t\rightarrow\pi{n}^2t$ substitutions in the Mellin transform corresponding to the $\Gamma$ function we have
\begin{equation}\label{RiCo}
\pi^{-s/2}\zeta(s)\Gamma(s/2)=\pi^{-s/2}\sum_{n=1}^\infty{n}^{-s}\int_0^\infty\frac{\d{t}}{t}t^{s/2}e^{-{t}}=\int_0^\infty\frac{\d{t}}{t}t^{s/2}\Psi(t)
\end{equation}
The most famous analytic continuation due to Riemann is by
${\Xi(s):=\pi^{-s/2}\Gamma(s/2)\zeta(s)}$ and
\begin{equation}\label{RiOr}\Xi\left(\frac{1+s}{2}\right)=-\frac{4}{1-s^2}+2\int_{1}^{\infty}\frac{\\\d{t}}{t}t^{1/4}\Psi\left(t\right)\cosh\left(\frac{\ln(t)s}{4}\right)\end{equation}
The \textit{\textbf{Poisson summation formula}} can be considered as a special case of a more geometric one called the Selberg trace formula and briefly is described as follows:

Let $F(x)$ for $x\rightarrow\pm\infty$ decays faster than $1/\vert{x}\vert^2$ (it is well-known that this assumption can be weakened in several ways). The series $\sum_{-\infty}^{\infty}F(x+n)$
is a well-defined function of the argument $x$ and has period $1$, hence can be expanded as a Fourier series and a calculation yields

\begin{equation}\label{PS1}
\sum_{n=-\infty}^{\infty}F(x+n)=\sum_{n=-\infty}^{\infty}\widehat{F}(n)e^{\I2\pi{n}{x}}
\end{equation}

The Theta $\Theta$ function is defined as the Poisson summation of a Gaussian, for instance we set
${\Theta({t}):=\sum_{n=-\infty}^{\infty}\exp(-\pi{n}^{2}{t})}$.
\ref{PS1} for $x=0$ combined with
$\widehat{\Di_{t}{f}}(p)=(\Di_{1/t}\widehat{f})(p)/\vert{t}\vert$
establishes the functional equation
$\Theta({t})=\frac{1}{\sqrt{{t}}}\Theta\left(\frac{1}{{t}}\right)$
because the Gaussian is a fixed point of the Fourier transform. 
This implies that $\Psi({t})=(\Theta({t})-1)/{2}$ enjoys 
${\Psi({t})=t^{-1/2}\left[\Psi\left({1}/{t}\right)+{1}/{2}\right]-{1}/{2}}$ and we finally split the integral on the r.h.s. of \ref{RiCo} at $1$ and elementary manipulations yield \ref{RiOr}.

In \ref{RiOr} it is obviously manifest that we have the celebrated functional equation $\Xi(s)=\Xi(1-s)$
Among other things \ref{RiOr} implies that $\zeta$ has so-called trivial zeros at the negative even integers $-2,-4,-6,-8,-10,\cdots$ because it is well-known that $\Gamma$ has simple poles at negative integers. Riemann conjectured \cite{RI} that all non trivial zeros of $\zeta$ have real part ${1}/{2}$.

\textbf{Acknowledgements:} I am deeply grateful to P. Moree for many inspiring discussions. I want to thank S. Bhattacharya, R. Friedrich, H. Furusho, D. Radchenko and O. Ramar\'e for reading draft versions, discussions and pointing out references.

\section{Riemann's illusive continuation formula}
\begin{definition}
Let $\alpha\in\mathbb{C}$ and the operators $H_\alpha,\Delta_\alpha:C^\infty({\mathbb{R}})\rightarrow{C}^\infty({\mathbb{R}})$ be defined by
\vspace{-0.05cm}$$H_\alpha:=\id+\alpha{t}\frac{\partial}{\partial{t}}\quad\text{and}\quad{\Delta}_\alpha:=2\alpha{t}\frac{\partial}{\partial{t}}+\alpha^2\left(t\frac{\partial}{\partial{t}}\right)^2$$
\end{definition}\vspace{-0.05cm}
For $\alpha=2$ we have essentially
${H}\propto{XP+PX}$, where $X$ and $P$ are the position and momentum operators respectively and the proportionality factor is $\I\hbar/2$ where $\hbar$ is the Planck constant. This Hamiltonian for the standard inner product was introduced by Connes in his approach towards a proof of the generalized RH see \cite{Co} and \cite{Co1} and the article \cite{BK} of Berry and Keating.

 A well-known fact is that $H_\alpha$ commutes with dilations $\Di_\beta$, {\em i.e.} 
$[H_{\alpha},\Di_{\beta}]=0$
because of the chain rule, hence $H_{\alpha'}$ commute with the Poisson summation maps $\Ps^\alpha$. For any $\alpha,\alpha'$ we can exchange $H_{\alpha}$ and $H_{\alpha'}$. By calculation 
we have $H_{\alpha}^2=\id+{\Delta}_\alpha$ and $H_\alpha$ also commutes with $\Delta_{\alpha'}$.

A crucial regularization property of an operator of the shape $H_\alpha=\id+\alpha{t}\partial_t$ is that it kills $t^{-1/\alpha}$ singularities. Hence a singularity of the shape $\sum_{i=1}^{k}c_i{t}^{-1/\alpha_i}$ gets canceled if we apply $\prod_{i=1}^{k}H_{\alpha_i}$. For a logarithmic singularity of the shape $t^{-1/\alpha}\ln^n(t)$ with $n\in\mathbb{N}$ we have
\vspace{-0.05cm}\begin{equation}\label{Help}
H^m_\alpha\left({t}^{-1/\alpha}\ln^{n}(t)\right)=m!\binom{n}{m}\alpha^m{t}^{-1/\alpha}\ln^{n-m}(t)
\end{equation}\vspace{-0.05cm}
 by a small calculation, hence we have the vanishing
$H^{n+1}_\alpha{t}^{-1/\alpha}\ln^n(t)=0$,
in other words $H_\alpha$ can also eliminate positive integer powers of logarithmic singularities.

This is a warning, the $\zeta$ function and Riemann's original continuation formula are quite intricate: For an unexperienced reader the l.h.s of Riemann's original illusive formula
\begin{equation}\label{FataMorgana1}
\frac{4}{1-z^2}=2\int_{1}^{\infty}\frac{\\\d{t}}{t}t^{1/4}\Psi(t)\cosh\left(\frac{\ln(t)z}{4}\right)
\end{equation}
for the $\Xi\left((1+z)/2\right)=0$ zeros $z$ and also the l.h.s of the formula
\begin{equation}\label{FataMorgana2}
\frac{4}{z-1/z}=2\int_{1}^{\infty}\frac{\\\d{t}}{t}t^{1/4}(H_4\Psi)(t)\sinh\left(\frac{\ln(t)z}{4}\right)
\end{equation}
may look like a direction to an algebraic coordinate descriptions of the $\Xi$ zeros, but this is more a {\em Fata Morgana} because in the further integration by parts procedure this hints get lost:

Let $\varphi:[1,\infty]\rightarrow\mathbb{C}$ be a smooth function, with no singularities at $1$ and in addition enjoying the property that at $\infty$ the function $\varphi$ and all of its derivatives are rapidly decaying.
\begin{proposition}\label{Step}
For $\varphi$ as above and $\pm1,0\neq{s}\in\mathbb{C}$ we have
\begin{align*}\int_{1}^{\infty}\frac{\d{t}}{t}t^{{1}/{\alpha}}\varphi\left(t\right)\begin{cases}\cosh\\\sinh\end{cases}\hspace{-0.4cm}\left(\frac{\ln(t)s}{\alpha}\right)&=\frac{-1}{s}\left[\begin{cases}0\\\alpha\varphi\left(1\right)\end{cases}\hspace{-0.4cm}+\int_{1}^{\infty}\frac{\d{t}}{t}t^{{1}/{\alpha}}(H_{\alpha}\varphi)(t)\begin{cases}\sinh\\\cosh\end{cases}\hspace{-0.4cm}\left(\frac{\ln(t)s}{\alpha}\right)\right]\end{align*}
\begin{align*}&\hspace{3.2cm}=\frac{-1}{1-s^2}\left[\begin{cases}(\alpha{H}_{\alpha}\varphi)(1)\\-s\alpha\varphi\left(1\right)\end{cases}\hspace{-0.4cm}+\int_{1}^{\infty}\frac{\d{t}}{t}t^{{1}/{\alpha}}({\Delta}_{\alpha}\varphi)(t)\begin{cases}\cosh\\\sinh\end{cases}\hspace{-0.4cm}\left(\frac{\ln(t)s}{\alpha}\right)\right]\end{align*}
\end{proposition}

We will use the full statement of the following auxiliary lemma \ref{Iteration} to show that it is not easy to chase the {\em Fata Morgana}, but here it is essentially enough to consider only the case $m=0$ and $\alpha=4$, for instance: If $\Psi:\mathbb{R}^{+}\rightarrow\mathbb{C}$ is smooth with
$\Psi(t)={t}^{-\frac{1}{2}}\left[\Psi\left(1/t\right)+{C}/{2}\right]-{C}/{2}$ the series
$\psi_{n}:=\Delta^n_{4}\Psi$ satisfies $(H_{4}\psi_0)(1)=-C/2$ and $(H_{4}\psi_n)(1)=0$ for all $n\in\mathbb{N}^+$ because
$\psi_n({1/t})={t}^{-\frac{1}{2}}\psi_n\left(t\right)+\delta_n^0{C}\bigr[{t}^{-\frac{1}{2}}-1\bigr]/{2}$
and $(H_{4}\psi_n)({1/t})=-{t}^{-\frac{1}{2}}(H_{4}\psi_n)\left(t\right)-C\delta_n^0\bigr[{t}^{-\frac{1}{2}}+1\bigr]/{2}$. 
\begin{proposition}\label{Iteration} Consider a smooth function $\Psi^\pm:\mathbb{R}^{+}\rightarrow\mathbb{C}$ enjoying
\begin{equation}\label{Funci1}
\Psi^\pm(t)=\pm(-1)^{m}{t}^{-\frac{2}{\alpha}}\Psi^\pm\left(1/t\right)+{C}\ln^{m}(t)\left[{t}^{-\frac{2}{\alpha}}\mp1\right]/{2}
\end{equation}
for some $m\in\mathbb{N}$ and $\alpha\in\mathbb{R}\setminus0$. The series $\psi_n$ defined by
$\psi^\pm_{n}:={\Delta}^n_{\alpha}\Psi^\pm$ respectively satisfy
\begin{align*}
&\psi^\pm_n(t)=\pm\frac{(-1)^m}{{t}^{\frac{2}{\alpha}}}\psi^\pm_n\left(\frac{1}{t}\right)+\frac{C\alpha^{2n}}{2}\sum_{i=0}^{n}\frac{(2n-i)!}{(\alpha/2)^i}\binom{m}{2n-i}\binom{n}{i}\ln^{m-2n+i}(t)\left[\frac{(-1)^i}{{t}^{\frac{2}{\alpha}}}\mp1\right]\\&
H_{\alpha}\psi^\pm_n(t)=\frac{\mp(-1)^{m}}{{t}^{\frac{2}{\alpha}}}H_{\alpha}\psi^\pm_n\left(\frac{1}{t}\right)\\&\quad\quad\quad\quad\quad-\frac{C\alpha^{2n}}{2}\sum_{i=0}^{n}\frac{(2n-i)!}{(\alpha/2)^i}\binom{m}{2n-i}\left(\binom{n+1}{i+1}+\binom{n}{i+1}\right)\ln^{m-2n+i}(t)\left[\frac{(-1)^i}{{t}^{\frac{2}{\alpha}}}\pm1\right]
\end{align*}
If $m$ is odd this implies $\psi^+_n(1)=-\frac{C{\alpha^{m}}}{2}{2^{2n-m}}m!\binom{n}{2n-m}
=-H_\alpha\psi^-_n(1)$
for all $n\in\mathbb{N}$ and if $m$ is even we have
${H}_\alpha\psi^+_n(1)=-\frac{C{\alpha^{m}}}{2}{2^{2n-m}}m!\left(\binom{n+1}{2n-m+1}+\binom{n}{2n-m+1}\right)=-\psi^-_n(1)$.
\end{proposition}
\begin{proof}[Proof] Notice \ref{Funci1} is self-consistent and $\psi(t)=\pm{t}^{-\frac{2}{\alpha}}\psi(1/t)$ implies $H_\alpha\psi(t)=\mp{t}^{-\frac{2}{\alpha}}(H_\alpha\psi)(1/t)$ by differential calculus. As a hint it is easier to evaluate $H_\alpha\Delta^n_\alpha{t}^{-\frac{2}{\alpha}}\ln^{m}(t)$ with help of \ref{Help} by switching operator coordinates, for instance put
$H_\alpha=2H_{\alpha/2}-\id$ and
$\Delta_\alpha=4\bigr(H_{\alpha/2}^2-H_{\alpha/2}\bigr)$
and use $(A+B)^n=\sum_{i=0}^{n}\binom{n}{i}A^iB^{n-i}$ for commuting operators. Now use $\binom{a+1}{b+1}=\binom{a}{b}+\binom{a}{b+1}$.\end{proof}

We consider \ref{Iteration} for the case $\alpha=4$, $n=0$: Riemann's analytic continuation \ref{RiOr} combined with the previous rewriting \ref{Step} shows the iteration of the integration by parts procedure is
\begin{equation}\label{IP}\Xi\left(\frac{1+s}{2}\right)=-\frac{4}{1-s^2}\delta^n_0+\left(\frac{-1}{1-s^2}\right)^n2\int_{1}^{\infty}\frac{\\\d{t}}{t}t^{1/4}({\Delta}_4^n\Psi)\left(t\right)\cosh\left(\frac{\ln(t)s}{4}\right)\end{equation}
$\forall{n}\in\mathbb{N}$. Formula \ref{IP} was discovered by P\'olya in the Nachlass of Jensen \cite{PO}.

Notice that
$({\Delta}_4\Psi)(t)\propto\sum_{n=1}^{\infty}\left[2\pi^2{n}^4{t}^2-3n^2{t}\right]\exp(-\pi{n}^{2}{t})$
is $[1,\infty]$ strictly positive, this contradicts zeros of $\Xi$ on $\mathbb{R}$. The function 
$({H}_4\Psi)(t)\propto\sum_{n=1}^{\infty}\left[1-4n^2{t}\right]\exp(-\pi{n}^{2}{t})$
is a strictly negative charged integral kernel restricted to the integration interval $[1,\infty]$.

Because of the functional equation
$({\Delta}^{n+1}_4\Psi)(t)=t^{-1/2}({\Delta}^{n+1}_4\Psi)(1/t)$
the functions $({\Delta}^{n+1}_4\Psi)(t)$ behave also nice for $t\rightarrow0$ and \ref{IP} implies
\begin{equation}\label{JAC}
\Xi\left(\frac{1+s}{2}\right)=\mathcal{M}\left[{\Delta}_4^{n+1}\Psi\right]\left(\frac{1+s}{2}\right)\Biggr/\big({s^2-1}\big)^{n+1}
\end{equation}
Formula \ref{IP} also implies that $\int_{1}^{\infty}{\\\d{t}}\;t^{-3/4}({\Delta}_4^n\Psi)\left(t\right)\cosh\left({\ln(t)s}/{4}\right)$ decays rapidly for $\vert\Im(s)\vert\rightarrow\infty$, roughly the oscillating in the integral kernel gets fast. The previous formula \ref{JAC} has been discussed by various people, and P\'olya's  conjecture that the Tur\'an inequalities are satisfied was verified \cite{Cs} by Csordas, Norfolk and Varga.

Let us follow \ref{FataMorgana1} approximately with help of \ref{Step}, \ref{Iteration}, the calculation is not important in the following but maybe a bit amusing: We first use the well-known hyperbolic addition theorem
$\cosh(a+b)=\cosh(a)\cosh(b)+\sinh(a)\sinh(b)$
to split up the real and imaginary part of the r.h.s. of the zero equation \ref{FataMorgana1}. For the imaginary part of \ref{FataMorgana1} we have with
$\sinh(x)=\sum_{m=0}^\infty{{x}^{2m+1}}/{(2m+1)!}$
the following formal expansion
$$\frac{8xy\I}{(1-x^2+y^2)^2+4x^2y^2}\stackrel{!}{=}\sum_{m=0}^\infty\frac{(x/4)^{2m+1}}{(2m+1)!}2\int_1^\infty\frac{\d{t}}{t}t^{\frac{1}{4}}\Psi(t)\ln^{2m+1}(t)\sinh\left(\frac{\ln(t)\I{y}}{4}\right)$$
\begin{align}\label{DeFata}
=\sum_{m=0}^\infty\frac{(x/4)^{2m+1}}{(2m+1)!}\Biggr[&4\I{y}\sum_{i=m+1}^{2m+1}\left(\frac{-1}{1+y^2}\right)^{i+1}4^{2m+1}2^{2i-2m-1}(2m+1)!\binom{i}{2m+1-i}\\&+\frac{2}{(1+y^2)^{2m+2}}\int_1^\infty\frac{\d{t}}{t}t^{\frac{1}{4}}\Delta_4^{2m+2}\bigr(\Psi(\cdot)\ln^{2m+1}(\cdot)\bigr)(t)\sinh\left(\frac{\ln(t)\I{y}}{4}\right)\Biggr]\nonumber
\end{align}
where we used \ref{Iteration} for the $+$ case and $\binom{a}{b}=\binom{a}{a-b}$. It would be absurd if we just truncate, approximate one side of \ref{DeFata} where the first line represents the only algebraic terms appearing and the second line integrals that decay faster than any negative power of $y$, because
\begin{equation}
\Delta_4^{2m+2}\bigr(\Psi(\cdot)\ln^{2m+1}(\cdot)\bigr)(t)=-t^{-1/2}\Delta_4^{2m+2}\bigr(\Psi(\cdot)\ln^{2m+1}(\cdot)\bigr)(1/t)
\end{equation}
The algebraic line can be rewritten as follows
$$\frac{-4\I{y}}{1+y^2}\sum_{m=0}^\infty{\left(\frac{x}{2}\right)^{2m+1}}\sum_{i=m+1}^{2m+1}\left(\frac{-4}{1+y^2}\right)^{i}\binom{i}{2m+1-i}$$
\begin{align}\label{DeFata}
&=\frac{-2\I{y}}{1+y^2}\sum_{i=0}^\infty\left(\frac{-2x}{1+y^2}\right)^{i}\sum_{m=0}^{\infty}{\left(\frac{x}{2}\right)^{m}}\binom{i}{m}\left(1-(-1)^{m+i}\right)\nonumber\\&=\frac{-2\I{y}}{1+y^2}\sum_{i=0}^\infty\left(\frac{-2x}{1+y^2}\right)^{i}\left[(x/2+1)^i-(-1)^i(-x/2+1)^i\right]\nonumber\\&=\frac{-2\I{y}}{1+y^2}\left[\frac{1}{1+\frac{x^2+2x}{1+y^2}}-\frac{1}{1+\frac{x^2-2x}{1+y^2}}\right]\nonumber=\frac{8xy\I}{(1-x^2+y^2)^2+4x^2y^2}\nonumber
\end{align}
and we have an agreement of the imaginary algebraic terms on the l.h.s and r.h.s of \ref{FataMorgana1}. Analog also the algebraic terms on both sides of the real part of illusion \ref{FataMorgana1} vanish self-consistent.

It may be also tempting for the reader to use the auxiliary lemmas \ref{Step} and \ref{Iteration} to rewrite parts of inequalities like the derivatives of $\Xi$ involving quite strange set of inequalities
\begin{equation}
0\leq\int_{1}^{\infty}\frac{\d{t}}{t}t^{1/4}f\bigr(\ln(t)\bigr)\begin{cases}\Psi\\-{H}_4\Psi\\\Delta_4\Psi\end{cases}\hspace{-0.4cm}(t)\begin{cases}\cosh\\\sinh\end{cases}\hspace{-0.4cm}\left(\frac{\ln(t)x}{4}\right)\left[\sqrt{2}\pm\cos\left(\frac{\ln(t)y}{4}\right)\pm'\sin\left(\frac{\ln(t)y}{4}\right)\right]\nonumber
\end{equation}
for all $f$ that are positive on $[0,\infty]$ and $x\geq0$ and substitute \ref{FataMorgana1} and \ref{FataMorgana2} in the inequalities at certain points of the calculation, but the author was not able to figure out something very meaningful on this way. However the reader is free to use his favourite trigonometric inequalities and try, but let us mention that the fact that $(H_4\Delta_4\Psi)(t)$ changes sign on $[1,\infty]$ is a crucial point for the proof of the subtle Tur\'an inequalities contained in \cite{Cs}.
\section{\textit{Mellin transforms and generic zeros}}\label{BEKE}
Formula \ref{JAC} can be considered as a special case of the following more general theorem \ref{PSC} corresponding to $H_2$ in the previous picture, in other words the following regularization is not the operator $H_4$ compatible with the symmetry ${\Psi({t})=t^{-1/2}\left[\Psi\left({1}/{t}\right)+{1}/{2}\right]-{1}/{2}}$ in the sense of \ref{Iteration}. Essentially \ref{PSC} maybe goes back to the year 1922, but to the best of our knowledge appears in the literature always in a slightly different form as we will explain and discuss briefly below the theorem:
\begin{theorem}\label{PSC}Let $0<\lambda\in\mathbb{R}$ and consider a function $f$ of rapid decay at $\infty$ and with $(H_1(f\circ\vert{\id}\vert^\lambda))\in\mathrm{C}^{2}(\mathbb{R}^+)$. We have an analytic continuation by the formula 
$$(1-s)\zeta(s)\mathcal{M}[f](s/\lambda)=\mathcal{M}\left[H_\lambda\circ{\Ps^\lambda}f\right](s/\lambda)$$
for $\Re(s)>0$, hence in this domain $(1-z)\zeta(z)=0\Rightarrow\mathcal{M}\left[H_\lambda\circ{\Ps^\lambda}f\right](z/\lambda)=0$.\end{theorem}\begin{proof}[Proof] We give a proof for the convenience of the reader, although the arguments are well-known: We define an even function $F_\lambda(t):=f(\vert{t}\vert^\lambda)$ and our assumptions on $f$ imply
\begin{equation}\label{PoissonSummation}
F_\lambda(0)+2\sum_{m=1}^\infty{F_\lambda\bigr(m\vert{t}\vert^{\frac{1}{\lambda}}\bigr)=f(0)+2\sum_{m=1}^\infty{f\bigr(m^\lambda{t}\bigr)}=\frac{1}{\vert{t}\vert^{\frac{1}{\lambda}}}\widehat{F_\lambda}(0)+\frac{2}{\vert{t}\vert^{\frac{1}{\lambda}}}\sum_{m=1}^\infty\widehat{F_\lambda}\left(\frac{m}{\vert{t}\vert^{\frac{1}{\lambda}}}\right)}
\end{equation}
or in other words we assumed that Poisson summation is valid.

It is well-known that the Fourier transform $\widehat{F}(p)$ satisfies $\widehat{F}(p)\leq{C}p^{-k}$ for real $p\rightarrow\infty$ with some constant $C$ if $g$ is $k$ times differentiable and if $F$ is smooth $\widehat{F}$ decays rapidly at $\infty$.

For a rapidly at $\infty$ decaying function $f(t)$ also $(\Ps^\lambda)f(t)=\sum_{n=1}^\infty{f}(n^{\lambda}{t})$ is rapidly decaying at $\infty$: A finite sum of rapidly decaying functions is of rapid decay and we have by assumption that for any $m\in\mathbb{N}$ exist $M_m,C_m\in\mathbb{N}$ with $\vert{f}(n^{\lambda}{t})\vert<C_{m}t^{-m}/n^{\lambda{m}}$ for $n\geq{M_m}$ hence if $m>1/\lambda$
$$\sum_{n=1}^\infty\vert{f}(n^\lambda{t})\vert\leq\sum_{n=1}^{M_m}\vert{f}(n^\lambda{t})\vert+C_{m}\zeta(\lambda{m}){t}^{-m}$$

We apply the rescaled operator $H_\lambda=\id+\lambda{t}\partial_t$ on the formula \ref{PoissonSummation} for $t\in\mathbb{R}^+$
\begin{equation}\label{PSID}
{F_\lambda(0)}+2\sum_{m=1}^\infty{(H_1{F_\lambda})\bigr(m{t}^{1/\lambda}\bigr)}=f(0)+2\sum_{m=1}^\infty{(H_\lambda{f})\bigr(m^\lambda{t}\bigr)}=\frac{2}{{t}^{1/\lambda}}\sum_{m=1}^\infty\bigr(\widehat{H_1{F_\lambda}}\bigr)\left(\frac{m}{{t}^{1/\lambda}}\right)
\end{equation}
where we used Leibniz rule, $H_\lambda{t}^{-1/\lambda}=0$, the fact that $t\partial_t$ is invariant under dilations, the formula
$t\partial_t(f\circ\id^\lambda)=\lambda(t\partial_t{f})\circ{\id}^\lambda$
and usual properties of the Fourier transform, for instance the identity $p\partial_p\widehat{F}(p)=-\widehat{\partial_x{x}F}(p)$.

Because $\lim_{t\rightarrow0}f(t)$ does not diverge by assumption we have for $\Re(s)>1$ the calculation
$$\zeta(s)\int_{0}^\infty\frac{\\\d{t}}{t}\;t^{s/\lambda}f(t)=\int_{0}^\infty\frac{\\\d{t}}{t}\;t^{s/\lambda}\sum_{n=1}^\infty{f({n^\lambda}t)}=-\frac{\lambda}{s}\int_{0}^\infty\frac{\\\d{t}}{t}\;t^{s/\lambda}t\partial_t({\Ps^\lambda}f)(t)$$
where we made the substitution $t={n^\lambda}t'$ and used $\sum_{n=1}^\infty{f({n^\lambda}t)}$ is of rapid decay at $\infty$.

Now consider the linear combination of integrals corresponding to $H_\lambda$: If $f\in\mathrm{C}^{3}(\mathbb{R}^+)$ the differentiation assumption in the theorem is obviously satisfied. Moreover we yield a continuation to $\Re(s)>0$ because under the assumptions $\widehat{H_1{F_\lambda}}(p)$ decays at least like $1/p^2$, hence
$$\frac{1}{\vert{t}\vert^{1/\lambda}}\sum_{m=1}^\infty\bigr(\widehat{H_1{F_\lambda}}\bigr)\left(\frac{m}{{t}^{1/\lambda}}\right)\leq{C}\zeta(2)\vert{t}\vert^{1/\lambda}$$
and this shows that the integral kernel on the r.h.s. of \ref{PSC} has no singularity as $t\rightarrow0$.\end{proof}
As mentioned the argument \ref{PSC} is essentially well-known since around 1922 and seems to go back to Eisenstein, Weil and M\"untz, but to the knowledge of the author appears in the literature generically in a slightly different, a bit more general shape: Let us again set $F_\lambda(t):=f(\vert{t}\vert^\lambda)$. The restrictions on the test functions in \ref{PSC} can be weakened in several ways, for example instead of introducing $\mathcal{H}$ we could continue to $\Re(s)>0$ with the weaker assumptions
\begin{equation}\label{AltCon}
F(0)=0\quad\text{and}\quad\widehat{F_\lambda}(0)=\int_{-\infty}^{\infty}\d{x}f(\vert{x}\vert^\lambda)=0
\end{equation}
hence the reader might wonder why we prefer a less general and simple statement with $H_\lambda\circ\Ps^\lambda$ and no integration assumptions. The first main advantage of this rewriting may be that we have an interpretation for the terms: $H_\lambda$ is a neat regularization of $\Ps^\lambda$, the composition is well-behaved.
Moreover as we will prospect in the following section \ref{hermit} our regularized combination is easy to symmetrize while the symmetrization of the condition \ref{AltCon} seems from our perspective somehow rigid and slightly less flexible considering standard symmetrization procedures, more precise the image under $\Ps^\lambda$ of functions satisfying the two conditions \ref{AltCon} does not necessarily satisfy \ref{AltCon}, only the first of this conditions, namely the vanishing at zero, is obviously conserved by $\Ps^\lambda$. Contrary the vanishing condition \ref{AltCon} is satisfied for the reminiscent Lie brackets
\begin{equation}
[f,g]_\lambda(t):=f(t)\int_{-\infty}^{\infty}\d{x}g(\vert{x}\vert^\lambda)-g(t)\int_{-\infty}^{\infty}\d{x}f(\vert{x}\vert^\lambda)
\end{equation}
Skew-symmetry of $[\cdot,\cdot]_\lambda$ is obvious and the validity of the Jacobi identity and $\int_{-\infty}^{\infty}\d{x}[f,g]_\lambda(\vert{x}\vert^\lambda)$ are just straight forward checks. It is also immediate that the space of functions satisfying \ref{AltCon} is an ideal with respect to $[\cdot,\cdot]_\lambda$. The condition  \ref{AltCon} is one of the main points where our description in subsection \ref{Std1} slightly differs from Connes spectral realization of the $\zeta$ zeros.

Let us also mention that if we suppose the stronger restriction that $f$ is smooth and of rapid decay at $\infty$ we can also continue to $\mathbb{C}$ quite analogous to Riemann's original calculation \ref{RiOr} by
\begin{equation}\label{AltCo}
\zeta(s)\int_{0}^\infty\frac{\d{t}}{t}\;t^{\frac{s}{\lambda}}f(t)=\frac{-\lambda}{2}\left[\frac{{F}_\lambda(0)}{s}+\frac{\widehat{F}_\lambda(0)}{1-s}\right]+\int_{1}^\infty\frac{\d{t}}{t}\sum_{n=1}^\infty\left[t^{\frac{s}{\lambda}}{F}_\lambda\bigr(n{t}^{\frac{1}{\lambda}}\bigr)+t^{\frac{1-s}{\lambda}}\widehat{F_\lambda}\bigr(n{t}^{\frac{1}{\lambda}}\bigr)\right]
\end{equation}

In consideration of the {\em Fata Morgana} \ref{FataMorgana1} and \ref{FataMorgana2} it seems interesting to try to  somehow reverse the integration by parts procedure \ref{Step} by inverting $H_{\alpha}^2=\id+{\Delta}_\alpha$ and technically cure the {\em a priori} quite ill substitution
$H_{\alpha}\stackrel{!}{=}\pm\sqrt{1+{\Delta}_{\alpha}}$.
Notice the Taylor expansion of $\sqrt{1+x}=\sum_{i=0}^\infty\binom{1/2}{i}x^i$ for $\vert{x}\vert<1$ only implies that this substitution is valid for test functions $t^\nu$ if
$\vert2\alpha\mu+\alpha^2\mu^2\vert<1$.

Notice that the operation sending $F(x)$ to $\sum_{-\infty}^{\infty}F(x+n)$ can be rewritten as the sum
$\id+\mathcal{R}+(-\id)\circ\mathcal{R}\circ(-id)$
where the operator $\mathcal{R}$ is defined by the formula
$(\mathcal{R}F)(x):=\sum_{n=1}^{\infty}F(x+n)$
and where we just set as usual
$((-id)f)(x)=f(-x)$. By this we want to capture that the map $\mathcal{R}$ defines a \textit{\textbf{Rota-Baxter operator}} of weight $1$, {\em i.e.} we have
\begin{equation}\label{RoBa}
\bigr(\mathcal{R}(F\cdot{G})\bigr)(x)=(\mathcal{R}F)(x)\cdot(\mathcal{R}G)(x)-\bigr(\mathcal{R}\left((\mathcal{R}F)\cdot{G}\right)\bigr)(x)-\bigr(\mathcal{R}\left(F\cdot(\mathcal{R}G)\right)\bigr)(x)
\end{equation}
Inspecting singularities this is paradox compared with \ref{PoissonSummation}: Poisson summation
$(\Ps^1f)(t)=(\mathcal{R}\circ\Di_t{f})(0)$
produces $\frac{1}{t}$ singularities on the l.h.s. while on the r.h.s. of \ref{RoBa} {\em a priori} a squared singularity of the shape $\frac{1}{t^2}$ appears and therefore should cancel with the other two more involved terms.
\subsection{P\'olya's extension of the Hermite-Biehler theorem}\label{PolyEx}
It seems interesting to try to extract out of \ref{PSC} information about the location of the $\zeta$ zeros by variation methods. The appearance of $\zeta$ in form $\Ps^\lambda$ is somehow responsible that there can be stable ``points": $\forall{s}\in\mathbb{C}$ the integral functional
$\varphi\rightarrow\int_{0}^\infty\frac{{\d}t}{t}\;t^{s/\lambda}(\varphi+\lambda{t}\varphi')$
does not admit stable ``points" $\varphi\in\mathrm{C}^{3}(\mathbb{R}^+)$, the Euler-Lagrange equations
$\frac{\d}{{\d}t}\frac{\partial{L}}{\partial{\varphi'}}-\frac{\partial{L}}{\partial\varphi}=0$
would imply
$(1-s)t^{s/\lambda-1}=0\;\forall{t}\;\in\mathbb{R}^+$
hence a contradiction if $s\neq1$. However, in the previous variation we neglect the more interesting non-local deformation $\Ps^\lambda$ representing the $\zeta$ function and just considered $H_\lambda$. Clearly $\Ps^\lambda$ is not a local operation we only have for the support of $H_\lambda\circ\Ps^\lambda{f}$ the implication
\begin{equation}\label{Support}
\mathrm{supp}(f)\subset[0,\delta]\Rightarrow\mathrm{supp}(H_\lambda\circ\Ps^\lambda{f})\subset[0,\delta]
\end{equation}

Obviously for some $\Lambda\in\mathbb{R}$ the symmetrized or skew-symmetrized integral functionals
$\varphi\rightarrow\int_{0}^\infty\frac{{\d}t}{t}\;t^{\Lambda/\lambda}(t^{s/\lambda}\pm{t}^{-s/\lambda})(H_\lambda\circ\Ps^\lambda\varphi)(t)$
respectively admit a stable vanishing imaginary or real part respectively if $\Re(s)=0$ and clearly a similar statement holds if $\Im(s)=0$.

Moreover P\'olya claims that the Euler product representation and the  \textit{\textbf{Hermite-Biehler theorem}} imply the observation
$\xi(1+s)+\xi(1-s)=0\Rightarrow{s}\in\I\mathbb{R}$
where
\begin{equation}\label{xi}
\xi\left(\frac{1+s}{2}\right):=(1-s^2)\Xi\left(\frac{1+s}{2}\right)=-2\int_1^\infty\frac{\d{t}}{t}t^{1/4}(\Delta_4\Psi)(t)\cosh\left(\frac{\ln(t){s}}{4}\right)
\end{equation}

In fact he proves a more general statement about how the zeros react if we impose quite the usual symmetric or also a skew-symmetric functional equation by brute force:
\begin{theorem}\label{PolYh}
If an entire real non-constant function $f(s)$ with Hadamard product
$$f(s)=c{e}^{-\gamma{s}^2+\beta{s}}z^q\prod\left(1-\frac{s}{z}\right){e}^{\frac{s}{z}}\left(1-\frac{s}{\overline{z}}\right){e}^{\frac{s}{\overline{z}}}$$
(where $c,\gamma,\beta$ are real constants, $\gamma\geq0$ and $q$ is an integer) has no zeros for $\Re(s)>\lambda$ then the two at the $\Lambda\geq\lambda$ axis respectively symmetrized or skew-symmetrized functions
$f(\Lambda+\I{s})\pm{f}(\Lambda-\I{s})$
respectively
vanish identically or admit only real zeros.
\end{theorem}
\begin{proof}[Proof]
Roughly the punch line of P\'olya's argumentation goes as follows, for details we refer to the original article \cite{PO}. First we can approximate $f$ by Hadamard factorization by polynomials, then the Hermite-Biehler theorem deduces the statement for the roots of the polynomials:

The Hermite-Biehler theorem states that a polynomial $p(x)=\mathrm{e}(x^2)+x\mathrm{o}(x^2)$ with real even part $\mathrm{e}(x^2)\in\mathbb{R}[x]$ and real odd part $\mathrm{o}(x^2)\in\mathbb{R}[x]$ is stable,{\em i.e.} $f(x)=0\Rightarrow\Re(x)<0$, if and only if $\mathrm{e}(-x^2)$ and $\mathrm{o}(-x^2)$ have simple real interlacing roots and $\Re\left(\mathrm{e}(x_0^2)/x_0\mathrm{o}(x_0^2)\right)>0$ for some $x_0$ with $\Re(x_0)>0$, it is well-known that this condition can also replaced by the condition that $[e,o]'(s)\geq0\forall{s}\in\I\mathbb{R}$ , here $[\cdot,\cdot]'$ is the so called Wronskian defined by $[f,g]':=f'g-fg'$, see \ref{DerLie}. Finally we conclude the theorem by Hurwitz theorem on the limit.\end{proof}

For example if $\xi$ has a zero free region $\Re(s)>\lambda$ for some $1/2\leq\lambda\in\mathbb{R}$ then $\forall\Lambda\geq\lambda$ we have
\begin{equation}
\xi(\Lambda+s)\pm\xi(\Lambda-s)\hspace{-0.1cm}=\hspace{-0.1cm}-4\hspace{-0.1cm}\int_1^\infty\hspace{-0.1cm}\frac{\d{t}}{t}t^{\frac{1}{4}}(\Delta_4\Psi)(t)\hspace{-0.1cm}\begin{cases}\cosh\\\sinh\end{cases}\hspace{-0.47cm}\left(\frac{\ln(t){(2\Lambda-1)}}{4}\right)\hspace{-0.1cm}\begin{cases}\cosh\\\sinh\end{cases}\hspace{-0.47cm}\left(\frac{\ln(t){s}}{2}\right)\hspace{-0.17cm}=0\hspace{-0.06cm}\Rightarrow{s}\in\I\mathbb{R}
\end{equation}
with the obvious exception that
$\xi\left(\frac{1}{2}+s\right)-\xi\left(\frac{1}{2}-s\right)$
identically vanishes, we refer to \cite{Lag} and \cite{Su1}.
\subsubsection{A comment on the Brujin-Newman-P\'olya operators}\label{Brujin}
In this context also other results of for example P\'olya \cite{PO}, \cite{PO1}, Brujin \cite{Bru}, Newman \cite{New} and Odlyzko \cite{Odl} on the reality of the roots of certain trigonometric integrals are interesting, we refer for more detailed considerations of the universal factors to\cite{KIMLEE}. Briefly they showed that certain classes of trigonometric integrals admit only real roots and that there exist universal integral kernel multiplication factors, in our notation  the functions $e^{-\rho\ln^2(t)}$
with $\rho\in\mathbb{R}^+$, conserving and even improving the reality of the roots of trigonometric integrals $\int_{0}^\infty\frac{{\d}t}{t}t^{\I{s}/2}\varphi(t)$ with the integral kernel symmetry
$\overline{\varphi(t)}=\varphi(1/t)$.

It is not difficult to project on solutions of this symmetry, namely for smooth $f$ of rapid decay at $0$ and $\infty$ the functions $f\rightarrow\varphi_f^\pm(t):=\big(f(t)\pm\overline{f(1/t)}\big)/2$ solve $\varphi^\pm(t)=\pm\overline{\varphi^\pm(1/t)}$ respectively. If $f$ is real and we have $f=\varphi_f^++\varphi_f^-$. The functions $e^{-\rho\ln^2(t)}$ are symmetric under $t\leftrightarrow1/t$ and it is again easy to project on the solutions, for instance for $f$ as above $f\rightarrow\phi_f^\pm(t):=\left(f(t)\pm{f}(1/t)\right)/2$ solves $\phi^\pm(t)=\pm\phi^\pm(1/t)$ and we have the decomposition $f=\phi_f^++\phi_f^-$.


Following P\'olya we have the Mellin transform interpretation
\begin{equation}\label{POLYI}
\mathcal{M}[e^{-\rho\ln^2(t)}](s)=\mathcal{M}[e^{-\rho\ln^2(t)}](0)\exp\left(\frac{s^2}{4\rho}\right)=\sqrt{\frac{\pi}{\rho}}\exp\left(\frac{s^2}{4\rho}\right)
\end{equation}
because we have by integration by parts the well-known differential equation
$$\partial_s\mathcal{M}[e^{-\rho\ln^2(t)}](s)=\int_{0}^\infty\frac{{\d}t}{t}t^{s}\ln(t)e^{-\rho\ln^2(t)}=\frac{-1}{2\rho}\int_{0}^\infty\hspace{-0.3cm}{{\d}t}\;t^{s}\partial_t{e}^{-\rho\ln^2(t)}=\frac{s}{2\rho}\mathcal{M}[e^{-\rho\ln^2(t)}](s)$$

Consider the functions
$\Xi_{f}(s):=\mathcal{M}[\Psi\cdot{f}](s/2)$ and $\tilde{\Xi}_{f}(s):=\mathcal{M}[(H_4\Psi)\cdot{f}](s/2)$
where we assume $f$ decays rapidly at $0$ and $\infty$ to avoid integration convergence issues, hence in this case $\Xi_f$ and $\tilde{\Xi}_f$ are obviously entire functions.

If $f$ is real we have $\overline{\Xi_{f}(s)}=\Xi_{f}(\overline{s})$. We can decompose $\Xi_{f}(s)=\Xi_{\phi_f^+}(s)+\Xi_{\phi_f^-}(s)$ and if $\phi^\pm(t)=\pm\phi^\pm(1/t)$ we can absorb the inhomogeneity $(t^{-1/2}-1)/2$ in the functional equation ${\Psi({t})=t^{-1/2}\left[\Psi\left({1}/{t}\right)+{1}/{2}\right]-{1}/{2}}$ in a quite symmetric way and for example rewrite $\Xi_{\phi^\pm}=\Lambda_{\phi^\pm}+\Omega_{\phi^\pm}$ with
$$\Lambda_{\phi^\pm}(s):=\int_1^\infty\frac{\d{t}}{t}\left(t^{\frac{s}{2}}\pm{t}^{\frac{1-s}{2}}\right)\Psi(t)\phi^\pm(t)\quad\text{and}\quad\Omega_{\phi^\pm}(s):=\int_0^1\frac{\d{t}}{t}\frac{t^{\frac{s-1}{2}}-t^{\frac{s}{2}}}{2}\phi^\pm(t)$$
clearly $\Lambda_{\phi^\pm}\hspace{-0.05cm}(s)\hspace{-0.05cm}=\hspace{-0.05cm}\pm\Lambda_{\phi^\pm}\hspace{-0.05cm}(1\hspace{-0.05cm}-\hspace{-0.05cm}s)$ but this symmetry is in general broken for the error summand $\Omega_{\phi^\pm}$:
\begin{proposition}\label{OddImp} Let  $\Xi^m_{\phi^\pm}(s):=\sum_{l=0}^m\Xi_{\phi^\pm}(s+l)$ and  $\tilde{\Xi}^m_{\phi^\pm}(s):=\sum_{l=0}^m(-1)^l\tilde{\Xi}_{\phi^\pm}(s+l)$.
\begin{align}\label{Para1}&\Xi^m_{\phi^\pm}(s)\mp\Xi^m_{\phi^\pm}(1-m-s)={\mathcal{M}\left[\frac{\phi^\pm}{2}\right]\hspace{-0.1cm}\left(\frac{s-1}{2}\right)-\mathcal{M}\left[\frac{\phi^\pm}{2}\right]\hspace{-0.1cm}\left(\frac{s+m}{2}\right)}\\&\tilde{\Xi}^m_{\phi^\pm}(s)\pm(-1)^m\tilde{\Xi}^m_{\phi^\pm}(s)(1-m-s)=-{\mathcal{M}\left[\frac{\phi^\pm}{2}\right]\hspace{-0.1cm}\left(\frac{s-1}{2}\right)-(-1)^m\mathcal{M}\left[\frac{\phi^\pm}{2}\right]\hspace{-0.1cm}\left(\frac{s+m}{2}\right)}\end{align}
\end{proposition}
\begin{proof}[Proof] The telescope identity $\Omega^m_{\phi^\pm}\left(s\right):=\sum_{l=0}^m\Omega_{\phi^\pm}(s+l)=\frac{1}{2}\int_0^1\frac{\d{t}}{t}\left(t^{\frac{s-1}{2}}-{t^{\frac{s+m}{2}}}\right){\phi^\pm}(t)$ holds. We also have the two parity equations $\mathcal{M}\left[{\phi^\pm}\right](s)=\pm\mathcal{M}\left[{\phi^\pm}\right](-s)$ respectively.
\end{proof}
The functions $\Xi^m_{e^{-\rho\ln^2(t)}}$ and $\tilde{\Xi}^m_{e^{-\rho\ln^2(t)}}$ solve $\left(\partial_\rho+4\partial^2_s\right)\Xi^m_{e^{-\rho\ln^2(t)}}(s)=0$, notice the similarity with the so-called heat equation. We have
$\Xi_{e^{-\rho\ln^2(t)}}(s)=\sum_{n=0}^\infty\frac{(-4\rho)^n}{n!}\partial^{2n}\Xi(s)$
and
$\tilde{\Xi}_{e^{-\rho\ln^2(t)}}(s)=(1-2\id)\Xi_{e^{-\rho\ln^2(t)}}(s)+16\rho\partial_s\Xi_{e^{-\rho\ln^2(t)}}(s)$ holds. \ref{OddImp} and \ref{POLYI} imply for $\forall{m}\in\mathbb{N}$ for example
$$\Omega^m_{e^{-\rho\ln^2(t)}}\left(s\right)-\Omega^m_{e^{-\rho\ln^2(t)}}\left({1-m-s}\right)={\mathcal{M}\big[e^{-\rho\ln^2(t)}\big](0)}\big(e^{{(s-1)^2}/{16\rho}}-e^{{(s+m)^2}/{16\rho}}\big)\hspace{-0.05cm}\big/2$$
and by this we yield for all $\rho>0$ the equivalences
\begin{align}\label{EquiSym1}
&\Xi^m_{e^{-\rho\ln^2(t)}}(s)-\Xi^m_{e^{-\rho\ln^2(t)}}(1-m-s)=0\Leftrightarrow{s}\in\frac{1-m}{2}+16\rho\pi\I\frac{\mathbb{Z}}{1+m}\\&{\tilde{\Xi}}^m_{e^{-\rho\ln^2(t)}}(s)+(-1)^m{\tilde{\Xi}}^m_{e^{-\rho\ln^2(t)}}(1-m-s)=0\Leftrightarrow{s}\in\frac{1-m}{2}+16\rho\pi\I\left(-\frac{1}{2}+\frac{\mathbb{Z}}{1+m}\right)\nonumber
\end{align}
Hence $\forall\rho>0$ solutions are on the line $\frac{1-m}{2}+\I\mathbb{R}$, if and only if $m$ is odd both symmetries on the l.h.s of \ref{EquiSym1} are satisfied simultaneous on the points
$\frac{1-m}{2}+16\rho\frac{\pi\I\mathbb{Z}}{1+m}$. Approximating by polynomials we also see that ${\partial^n_s\Xi}^m_{e^{-\rho\ln^2(t)}}(s)-(-1)^{n}{\partial^n_s\Xi}^m_{e^{-\rho\ln^2(t)}}(1-m-s)=0\Rightarrow{\Re(s)=1/2}$ by the Gauss-Lucas theorem. Unfortunately $\lim_{\rho\rightarrow0}\Xi_{e^{-\rho\ln^2(t)}}=\Xi$ and in the limit \ref{EquiSym1} vanishes because of $\Xi(1-s)=\Xi(s)$. The previous considerations seem a bit bizarre in consideration of \ref{PolYh} and the remarkable result of Speiser that a zero of $\zeta$ is equivalent to a zero of the derivative $\zeta'$ inside the critical strip \cite{Sp}.

In the sequel to the present paper we will discuss the more general family of integrals $$\Xi(\rho,\vec{s}\;):=\int_{0}^\infty\frac{{\d}t_1}{t_1}\cdots\int_{0}^\infty\frac{{\d}t_n}{t_n}\left(\prod_{i=1}^n{t^{{s_i}/2}_i}\Psi(t_i)
\right)e^{-\sum_{1\leq{i,j}\leq{n}}\rho_{ij}\ln(t_i)\ln(t_j)}$$
where $\rho$ denotes a symmetric $n\times{n}$ matrix underlying certain convergence restrictions and $\vec{s}\in\mathbb{C}^n$.
\section{\textit{Some comments on the Hilbert-P\'olya conjecture}}\label{hermit}
The \textit{\textbf{Hilbert-P\'olya conjecture}} goes back to a discussion of P\'olya with Landau and it states that the $\Xi$ zeros should be interpreted as the energy levels of a quantum system. More precise the conjecture is the existence of a hermitian operator $\textit{\textbf{H}}$ with the property that if
$\Xi(1/2+\I{E})=0$
then there is a eigenstate $\psi_E$ with eigenvalue $E$, {\em i.e.}
$\textit{\textbf{H}}\psi_E=E\psi_E$. Because the absolute value of non-trivial zeros can be arbitrary big by \cite{PO} the Hilbert-P\'olya operator $\textit{\textbf{H}}$ can only be realized as an unbounded operator. Works of Berry, Keating and Connes gave various evidence that $\textit{\textbf{H}}$ should be some quantization of $H=XP+PX$. One of the ideas in \cite{BK} is to get a discrete spectrum by imposing boundary conditions for the space of functions on which $H$ acts. A keystone in \cite{Co} is to incorporate a topology related to prime numbers and from this point of view our basic considerations might be naive.

It is well-known that the Schwartz space is a dense subspace of the Hilbert space $L^2(\mathbb{R})$ of square integrable functions. The completion of pre-Hilbert spaces and the extension of symmetric operators are an intensive studied research topic where the Cayley transform, {\em i.e.} the map  $A\rightarrow\mathrm{U}_A:=(A-\I)\circ(A+\I)^{-1}$ taking a symmetric operator $A$ to the isometry $\mathrm{U}_A:\mathrm{Ran}(A+\I)\rightarrow\mathrm{Ran}(A-\I)$, provides useful reformulations. In the following we work with an adequate not Cauchy complete space $\mathcal{H}$ equipped with an inner product:
\begin{definition}\label{SP} Consider the pre-Hilbert space $\mathcal{H}$ of rapidly at $0$ and $\infty$ decaying smooth functions equipped with the inner product
${\langle{f_1},f_2\rangle_\Def:=\int_0^\infty\d{t}\;t^\Def
\overline{f_1}(t)f_2(t)}$
for $\Def\in\mathbb{R}$.\end{definition}

Consider a function $V$ that is in addition real on the real axis, {\em i.e.} $V(t)=[v(t)+\overline{v(t)}]/2=\overline{V(t)}$ holds $\forall{t}\in\mathbb{R}^+$, it is well-known that for such a real $V$ the operator $\Mu_V$ is hermitian where
\begin{definition}\label{MuLti} For a function $V$ with at worst $t^\lambda$ singularities at $0$ and $\infty$ let as usual the multiplication operator $\Mu_V:\mathcal{H}\rightarrow\mathcal{H}$ operate on states $f\in\mathcal{H}$ by ${f(t)\rightarrow\left(\Mu_V{f}\right)(t):={V}(t)\cdot{f}(t)}$.
\end{definition}
\begin{definition}
Let $\mathfrak{g}:=\{A:\mathcal{H}\rightarrow\mathcal{H}\vert\exists{A}^\ast:\mathcal{H}\rightarrow\mathcal{H}:\langle{f},A{g}\rangle_\Def=\langle{A^\ast{f}},{g}\rangle_\Def\;\forall{f,g}\in\mathcal{H}\}$ denote the (Lie-)algebra of operators admitting an adjoint and $\mathfrak{s}$ the sub Lie-algebra of Hamiltonians.
\end{definition}

$\Ps^\lambda$ begins with the identity, by the previous reasonings is related to $\zeta$ and in \ref{Ad-Hoc} we introduce a non-local deformation of $H_\lambda$ with help of $\Ps^\lambda$, for the cook we just need one more ingredient:
\begin{definition}\label{defiota}
By ${\tau_\Def}$ we denote $\tau_\Def^{-1}$ where for $\mu\in\mathbb{R}\setminus0$ the operators $\tau_\Def^\mu$ are defined by
$$({\tau_\Def^\mu}{f})(t):=\bigr((\Mu_{t^{1+\Def}}{f})\circ{\id}^{\mu}\bigr)(t)=t^{\mu(1+\Def)}f\bigr(t^{\mu}\bigr)$$
\end{definition}

Notice that for $\Def=0$ we have defined the involution appearing in Weil's positivity criteria \ref{Weil} and that the previous considerations \ref{Iteration} and \ref{Brujin} can shortly be rewritten with this definition, we have for example the formula
$({\tau_\Def^\mu}\circ{\tau_\Def}{f})(t)=(f\circ\id^{-\mu})(t)=(\tau^{-\mu}_{-1}f)(t)=f(1/t^\mu)$. It is clear that $t^{\frac{\mu(1+\Def)}{1-\mu}}$ is a fixed point of ${\tau_\Def^\mu}$ and some of the following reasonings in summary show that ${\tau_\Def}$ for the test functions $\mathcal{H}$ is a hermitian involution, {\em i.e.}
${\tau_\Def}^\ast={\tau_\Def}$ and ${\tau_\Def}\circ{\tau_\Def}=\id$. 

In \ref{defiota} we have not defined an honest group action of $\mathbb{R}\setminus0$, the composition just satisfies
\begin{equation}\label{mult}
\bigr({\tau_\Def^{\mu_1}}\circ{\tau_\Def^{\mu_2}}{f}\bigr)(t)=\bigr(\Mu_{t^{\mu_1(1+\Def)}}\circ{\tau_\Def}^{\mu_1\mu_2}{f}\bigr)(t)=\Bigr({\tau_\Def^{\mu_1\mu_2}}\circ\Mu_{t^{\frac{1+\Def}{\mu_2}}}{f}\Bigr)(t)=t^{\mu_1(1+\Def)}\bigr({\tau_\Def^{\mu_1\mu_2}}f\bigr)(t)
\end{equation}
and hence only a bona fide group action of $\mathbb{R}\setminus0$ if $\Def=-1$. The inverse $({\tau_\Def^\mu})^{-1}$ of ${\tau_\Def^\mu}$ is given by
$\bigr(({\tau_\Def^\mu})^{-1}f\bigr)(t):=\Bigr(\Mu_{t^{-\frac{(1+\Def)(1+\mu)}{\mu}}}\circ{\tau_\Def^{1/{\mu}}}{f}\Bigr)(t)=\frac{1}{t^{1+\Def}}f\bigr(t^{1/\mu}\bigr)$.
A reason for the definition \ref{defiota} can be seen in the adjunction 
$\langle{\tau_\Def^\mu}{f_1},f_2\rangle_\Def=\langle{f_1},({\tau_\Def^{\mu}})^\ast{f_2}\rangle_\Def$
with the adjoint $({\tau_\Def^{\mu}})^\ast$ given by
$\bigr({\tau_\Def^\mu}\bigr)^\ast={{\tau_\Def^{1/\mu}}}/{\vert\mu\vert}$
where we restrict to test functions in $\mathcal{H}$ and with respect to \ref{SP}.

Let the adjoint conjugation be defined by
$\mathrm{Con}^\ast_{B}(A):=B\circ{A}\circ{B}^\ast$. It is clear by essentially iteration of the adjoint conjugation
$A\rightarrow{B}\circ{A}\circ{B}^\ast$
of $A\in\mathfrak{s}$ with $B\in\mathfrak{g}$ that for a sequence $\mu_1,\cdots,\mu_n\in\mathbb{R}\setminus0$ we have
${\tau_\Def}^{\mu_1}\circ\cdots\circ{\tau_\Def}^{\mu_n}\circ{A}\circ{\tau_\Def}^{1/\mu_n}\circ\cdots\circ{\tau_\Def}^{1/\mu_1}\in\mathfrak{s}$
but we still did not define by this an honest action of $\mathbb{R}\setminus0$ on $\mathfrak{s}$: Explicit we just have by \ref{mult}
\begin{align}\label{pseudogroup}
{\tau_\Def}^{\mu_1}\circ{\tau_\Def}^{\mu_2}\circ{A}\circ{\tau_\Def}^{1/\mu_2}\circ{\tau_\Def}^{1/\mu_1}&={\tau_\Def}^{\mu_1\mu_2}\circ\Mu_{t^{\frac{1+\Def}{\mu_2}}}\circ{A}\circ\Mu_{t^{\frac{(1+\Def)}{\mu_2}}}\circ{\tau_\Def}^{\frac{1}{\mu_1\mu_2}}\\&=\Mu_{t^{\mu_1(1+\Def)}}\circ{\tau_\Def}^{\mu_1\mu_2}\circ{A}\circ{\tau_\Def}^{\frac{1}{\mu_1\mu_2}}\circ\Mu_{t^{\mu_1(1+\Def)}}\nonumber
\end{align}
and hence again only an action of the multiplicative group $\mathbb{R}\setminus0$ for $\Def=-1$. 

The involution $\tau_\Def$ also admits n-th roots, for example $\forall{k}\in\mathbb{Z}$ the n-th power of the maps $\frac{e^{\I\pi(1+2k)/n}-1}{e^{\I\pi(1+2k)/n}+1}+\tau_\Def$ are proportional to $\tau_\Def$, in particular the maps $(\I\pm\tau_\Def)/\sqrt{\pm\I2}$ square to $\tau_\Def$.

The operators $(\id\pm{\tau_\Def})$ satisfy $(\id+{\tau_\Def})\circ(\id-{\tau_\Def})=0$, $(\id\pm{\tau_\Def})^m=2^m(\id\pm{\tau_\Def})$ and $(\id\pm{\tau_\Def})\circ{\tau_\Def}=\pm(\id\pm{\tau_\Def})$. We can split $f\in\mathcal{H}=\mathcal{H}_\Def^+\oplus\mathcal{H}_\Def^-$ in
$f=f^{+}+f^{-}$
where we project
\begin{equation}\label{IotaProjection}
f\rightarrow{f}^{\pm}:=(f\pm{\tau_\Def}{f})/2\in\mathcal{H}^\pm_\Def
\end{equation}
onto the two eigenspaces of ${\tau_\Def}$, hence ${\tau_\Def}{f}^{\pm}=\pm{f}^{\pm}$ is valid and $\langle{f^+},{f^-}\rangle_\Def=0$.

The operator $\Mu_{\ln(t)}$ satisfies $\Mu_{\ln(t)}\circ\tau_\Def=-\tau_\Def\circ\Mu_{\ln(t)}$ and hence defines maps $\mathcal{H}_\Def^+\leftrightarrow\mathcal{H}_\Def^-$, but admits no inverse because $\Mu_{1/\ln(t)}$ could produce a singularity at $1$ and is not in any case well-defined, $\Mu_{\ln(t)}$ maps $\mathcal{H}_\Def^-$ to functions in $\mathcal{H}_\Def^+$ that in addition vanish at $1$. We have $[t\partial_t,\Mu_{\ln(t)}]_{-}=\id$.

Let us refer to section \ref{Hdependence} for a more detailed discussion of the $\Def$-dependence and just mention here that the two projections in \ref{IotaProjection} transform according to the rule
$\id\pm\tau_{\Def'}={\mathrm{Con}_{\Mu_{t^{\frac{\Def-\Def'}{2}}}}}(\id\pm\tau_{\Def})$.

With help of $\langle{\tau_\Def^\mu}{f_1},f_2\rangle_\Def=\langle{f_1},({\tau_\Def}^{\mu})^\ast{f_2}\rangle_\Def$ by elementary manipulations we have the formula
\begin{equation}\label{IotaM}
\mathcal{M}\bigr[\tau_\Def^{\mu}{f}\bigr](s)=\frac{1}{\vert\mu\vert}\mathcal{M}[f]\left(1+\Def+\frac{s}{\mu}\right)\end{equation}
Let $\;\widehat{}\;$ again denote the usual Fourier-transform. On $\mathcal{H}$ the regularized Poisson summation \ref{PSID} can be rewritten with help of $\tau$ in the shape
$H_1\circ\Ps^1=\tau^{-1}_0\circ\Ps^1\circ\widehat{}\circ{H_1}$
and we have $H_\lambda\circ\Ps^\lambda=\tau^{1/\lambda}_{-1}\circ{H}_1\circ\Ps^1\circ\tau^{\lambda}_{-1}$.

\subsection{\textit{Proof I, the four standard symmetrizations of $H_\lambda\circ{\Ps^\lambda}$}}\label{Std1}
\vspace{-0.1cm}\begin{definition}\label{SPpm} Consider the pre-Hilbert space $\mathcal{H}_{\infty,0}$ of at $\infty$ or at $0$ rapidly decaying smooth functions equipped with the inner product
$\langle{f_1},f_2\rangle_\Def:=\int_0^\infty{\d{t}}\;t^\Def\overline{f_1}(t)f_2(t)$
for $\Def\in\mathbb{R}^\pm$ respectively.\end{definition}
In consideration of ${\tau_\Def}$ there is some tension, we have ${\tau_\Def}\mathcal{H}_\infty\subset\mathcal{H}_0$ but a function in ${\tau_\Def}\mathcal{H}_0$ could have a $t^{-(1+\Def)}$ singularity at zero, and we want to consult this by \ref{SPpm}.

Adjunction $\;^\ast=\;^{\ast(\Def)}$ is an involution and we invoke a natural projector $\mathrm{p}=\mathrm{p}(\Def)$ and square zero map $\partial=\partial(\Def)=\mathrm{p}\circ\Mu_{\I}$ where we suppress in our notation the $\Def$ dependence:
\begin{align}\label{TriDi}
&A\in\mathfrak{g}\rightarrow\mathrm{p}(A):=({A}+A^{\ast})/2\in\mathfrak{s}\quad\text{and}\quad{A}\in\mathfrak{g}\rightarrow\partial({A}):=\I({A}-A^{\ast})/2\in\mathfrak{s}
\end{align}
We can translate between $\langle\cdot,\cdot\rangle_\Def$ and $\mathcal{M}$ by
$\langle{f_1},\Mu_{t^{\I{y}}}f_2\rangle_\Def=\mathcal{M}[\overline{f}_1\cdot{f}_2](\Def+1+\I{y})$ and have
\begin{equation}
\left\langle{f_1},\begin{cases}{\mathrm{p}}\\{\partial}\end{cases}\hspace{-0.35cm}\Mu_{t^{\I{y}}}\;f_2\right\rangle_\Def=\begin{cases}{\Re}\\{\Im}\end{cases}\hspace{-0.4cm}\Bigr(\mathcal{M}[\overline{f}_1\cdot{f}_2](\Def+1+\I{y})\Bigr)
\end{equation}
Notice $\id=\mathrm{p}-\I\partial$ and$\;^\ast=\mathrm{p}+\I\partial$ and the restriction identities $\mathrm{p}\vert_{\mathfrak{s}}=\id$, $\partial\vert_{\mathfrak{s}}=0$ hold, hence $\mathrm{p}^2=\mathrm{p}$, $\partial^2=\partial\circ\mathrm{p}=0$ and $\mathrm{p}\circ\partial=\partial$. To call $\partial$ a differential is misleading, we just have $\partial[A,{B}]_{-}=\frac{\I}{2}[\mathrm{p}A,\mathrm{p}B]_{-}-\frac{\I}{2}[\partial{A},\partial{B}]_{-}$ and $\mathrm{p}[A,{B}]_{+}=\frac{1}{2}[\mathrm{p}A,\mathrm{p}B]_{+}-\frac{1}{2}[\partial{A},\partial{B}]_{+}$. and clearly $\mathrm{ker}\partial/\mathrm{im}\partial=0$.
\begin{theorem}\label{Ad-Hoc}
Consider the pre-Hilbert space $\mathcal{H}$ equipped with the inner product $\langle\cdot,\cdot\rangle_\Def$, see \ref{SP}. For $0<\lambda\in\mathbb{R}$ the operators
$\mathcal{Z}^\lambda_\pm:\mathcal{H}\rightarrow\mathcal{H}$
defined by
\begin{equation}\label{DH}
\mathcal{Z}^\lambda_\pm:=\begin{cases}\mathrm{p}\left({H}_\lambda\circ\Ps^\lambda\right)=\bigr(H_\lambda\circ{\Ps^\lambda}+{\tau_\Def}\circ{H_\lambda}\circ{\Ps^\lambda}\circ{\tau_\Def}\bigr)/2\\\partial\left({H}_\lambda\circ\Ps^\lambda\right)=\I\bigr(H_\lambda\circ{\Ps^\lambda}-{\tau_\Def}\circ{H_\lambda}\circ{\Ps^\lambda}\circ{\tau_\Def}\bigr)/2\end{cases}
\end{equation}
are self-adjoint. We have $[\mathcal{Z}^\lambda_\pm,\tau_\Def]_\mp=0$ and hence restrictions $\mathcal{Z}^\lambda_\pm:\mathcal{H}_\Def^{\pm'}\rightarrow\mathcal{H}_\Def^{\pm'\pm}$ with respect to the two eigenspaces of $\tau_\Def$. Restricted to $\mathcal{H}$ we also have hermitian operators by the formulas
\begin{equation}\label{DH2}
\widehat{\mathcal{Z}}^\lambda_\pm:=\begin{cases}{H}_\lambda\circ\Ps^\lambda\circ{\tau_\Def}:{\tau_\Def}\mathcal{H}_\infty\rightarrow\mathcal{H}_\infty\\{\tau_\Def}\circ{H_\lambda}\circ\Ps^\lambda:\mathcal{H}^\infty\rightarrow{\tau_\Def}\mathcal{H}_\infty\end{cases}
\end{equation}\end{theorem}
\begin{proof}[Proof] For $H_\alpha$ the adjoint with respect to the inner product \ref{SP} is given by 
\begin{equation}\label{Hadjoint}
{H}_\alpha^\ast=\bigr(1-\overline{\alpha}(1+\Def)\bigr){H}_{\frac{\overline{\alpha}}{\overline{\alpha}(1+\Def)-1}}\end{equation}
Hence for $\Def\neq-1$ with respect to the inner product $\langle\cdot,\cdot\rangle_\Def$ the operator
${H}_{\frac{2}{1+\Def}}/2\I$
is self-adjoint and for the special value $\Def=-1$ the operator
$\I{t}\partial_t$
is self-adjoint resolving the for $\Def=-1$ singular appearing expression ${H}_{\frac{2}{1+\Def}}/2\I$ consistently. It is standard that for $A\in\mathfrak{g}$ with adjoint $A^\ast\in\mathfrak{g}$ we have $A^\ast\circ{A},A\circ{A}^\ast\in\mathfrak{s}$ and also real polynomials of a self-adjoint operator $A\in\mathfrak{s}$ are self-adjoint. Hence for $\Def\neq-1$ we have
$\Delta_{\frac{2}{1+\Def}}\in\mathfrak{s}$
and for the special value $\Def=-1$ we have
$({t}\partial_t)^2\in\mathfrak{s}$.

With respect to \ref{SP} we have for $\Di_\beta$ the adjoint operator by the enlightening rewriting
\begin{equation}\label{Dadjoint}
\Di_\beta^\ast=\frac{1}{\beta^{1+\Def}}\Di_{1/\beta}={\tau_\Def}\circ\Di_{\beta}\circ{\tau_\Def}
\end{equation}
The compositions of operators ${\tau_\Def^\mu}$ and $\Di_\beta$ transform by
${\tau_\Def^\mu}\circ\Di_{\beta}=\frac{1}{\beta^{1+\Def}}\Di_{\beta^{1/\mu}}\circ{\tau_\Def}^{\mu}$ and we find
${\tau_\Def^\mu}\circ\Ps^\alpha=\left(\sum_{n=1}^\infty\frac{1}{n^{\alpha(1+\Def)}}\Di_{n^{\alpha/\mu}}\right)\circ{\tau_\Def^\mu}$.
Observe $\Ps^\alpha$ is not hermitian: We have
$\langle\Ps^\alpha{f}_1,f_2\rangle_\Def=\langle{f}_1,(\Ps^\alpha)^{\ast}{f}_2\rangle_\Def$
with the twist
\begin{equation}
(\Ps^\alpha)^{\ast}=\sum_{n=1}^\infty\frac{1}{n^{2}}\Di_{n^{-\alpha}}={\tau_\Def}\circ\Ps^\alpha\circ{\tau_\Def}
\end{equation}
Between the operators ${\tau_\Def^\mu}$ and $H_\alpha$ by calculation the commutation relations
\begin{align}\label{ComRel1}
&{{H}_{\alpha}\circ{\tau_\Def^{\mu}}=\bigr(1+\alpha\mu(1+\Def)\bigr){\tau_\Def^{\mu}}\circ{H}_{\frac{\alpha\mu}{1+\alpha\mu(1+\Def)}}}\\&{{\tau_\Def^{\mu}}\circ{H}_{\alpha}=\bigr({1-\alpha(1+\Def)}\bigr){H}_{\frac{\alpha}{\mu\left(1-\alpha(1+\Def)\right)}}}\circ{\tau_\Def^{\mu}}\nonumber
\end{align}
are valid. The formulas \ref{ComRel1} imply that we have 
$\tau_\Def^{\frac{-2/\mu}{1-2\Def}}\circ{H}_{\frac{2}{1+\Def}}=\frac{\Def-1}{\Def+1}{H}_{\frac{2}{1+\Def}}\circ{\tau_\Def^{\frac{-2}{\mu(1+\Def)}}}$, ${H}_{\frac{2}{1+\Def}}\circ{\tau_\Def^{\mu}}=(1+2\mu){\tau_\Def^{\mu}}\circ{H}_{\frac{2\mu}{(1+\Def)(1+2\mu)}}$ and ${H}^m_{\frac{\mu-1}{\mu(1+\Def)}}\circ\left({\tau_\Def^{\mu}}\right)^n=\mu^{m\cdot{n}}\left({\tau_\Def^{\mu}}\right)^{n}\circ{H}^m_{\frac{\mu-1}{\mu(1+\Def)}}$ where we suppose $\mu\neq0$, $\Def\neq-1$ respectively. This shows in particular the anti-commuting
$H_{2/(1+\Def)}\circ{\tau_\Def}=-{\tau_\Def}\circ{H}_{2/(1+\Def)}$ and hence the rewriting $\Delta_\alpha=(\I\tau_{\frac{2-\alpha}{\alpha}}\pm{H}_\alpha)^2$.

With the preparation formulas $\bigr({\tau_\Def^\mu}\bigr)^\ast={{\tau_\Def^{1/\mu}}}/{\vert\mu\vert}$, \ref{Dadjoint}, \ref{Hadjoint} and \ref{ComRel1} we have for $\lambda\in\mathbb{R}$
\begin{equation}\label{hermitri}
\bigr((H_\lambda\circ\Ps^\lambda)^n\bigr)^\ast=\bigr((H_\lambda\circ\Ps^\lambda)^\ast\bigr)^n={\tau_\Def}\circ(H_\lambda\circ\Ps^\lambda)^n\circ{\tau_\Def}
\end{equation}
 by induction where we essentially only use $(A\circ{B})^\ast=B^\ast\circ{A}^\ast$, \ref{Dadjoint}, ${\tau_\Def}^\ast={\tau_\Def}$, ${\tau_\Def}\circ{\tau_\Def}=\id$ and
$$\frac{\id}{(1+\Def)\id-1}\circ\frac{\overline{\id}}{(1+\Def)\overline{\id}-1}=\overline{\id}$$
for the composition of the respective M\"obius-transforms. Hence the two standard procedures \ref{TriDi} lead for all $<\lambda\in\mathbb{R}$ to the formulas \ref{DH} where it is manifest that the two operators $\mathcal{Z}^\lambda_\pm$ are well-defined maps $\mathcal{H}\rightarrow\mathcal{H}$, compare the argumentation in the proof of \ref{PSC}. It is clear that by \ref{hermitri} the standard symmetrizations $\mathrm{p}$ and $\partial$ are compatible with the dynamics generated by $H_\lambda\circ\Ps^\lambda$, {\em i.e.} the iterations $(H_\lambda\circ\Ps^\lambda)^n$ behave similar under the symmetrizations $\mathrm{p}$ and $\partial$, the same holds for the other two symmetrization procedures \ref{DH2}: The fact that $\widehat{\mathcal{Z}}^\lambda_\pm$ are hermitian can be justified in two ways, the first method is just a direct check by calculation. We have 
$$({\tau_\Def}\circ({H}_\lambda\circ\Ps^\lambda)^n)^\ast=\left({\tau_\Def}\circ(\Ps^\lambda\circ{H}_\lambda)^n\circ{\tau_\Def}\right)\circ{\tau_\Def}={\tau_\Def}\circ({H}_\lambda\circ\Ps^\lambda)^n$$
because of \ref{hermitri} and ${\tau_\Def}^\ast={\tau_\Def}$. The adjoint conjugation $\mathrm{Con}^\ast_{B}(A)=B\circ{A}\circ{B}^\ast$ is also self-adjoint
if $A\in\mathfrak{s}$, hence
${H}_\lambda\circ\Ps^\lambda\circ{\tau_\Def}=\mathrm{Con}_{\tau_\Def}({\tau_\Def}\circ{H}_\lambda\circ\Ps^\lambda)\in\mathfrak{s}$
 if ${\tau_\Def}\circ{H}_\lambda\circ\Ps^\lambda\in\mathfrak{s}$ because ${\tau_\Def}\circ{\tau_\Def}=\id$. The second somehow slightly more conceptional proof will be discussed in section \ref{proof2}.\end{proof}
A strategy developed in \cite{ConMar},\cite{RM} is to consider the dual of the quotient of a certain larger function space by the range of $\Ps^\lambda$ and the transposed operators acting on the dual space. This procedure yields a spectral realization of the Riemann zeros, in analogy we will consider the quotients by the image of the operators \ref{Ad-Hoc} in the following section \ref{QuoEigen}. Because in \ref{EigenState1},\ref{EigenState2} we do not consider the mentioned quotient and the transposition of $t\partial_t$ we have still an inner product but also corrections to the usual spectral realization.

\subsubsection{\textit{Eigenstates of $\mathcal{Z}^\lambda_\pm$ and some sums of $\zeta$ functions}}\label{EigenState1}
In the section \ref{PSC} we showed
$\mathcal{M}[H_\lambda\circ{\Ps^\lambda}f]\left(\frac{s}{\lambda}\right)=(1-\id)\zeta(s)\mathcal{M}[f]\left(\frac{s}{\lambda}\right)$
and \ref{IotaM} applied twice yields $\mathcal{M}[{\tau_\Def}\circ{H_\lambda}\circ{\Ps^\lambda}\circ{\tau_\Def}{f}]\left(\frac{s}{\lambda}\right)=(1-\id)\zeta\left((1+\Def)\lambda-s\right)\mathcal{M}[{f}]\left(\frac{s}{\lambda}\right)$, hence
by linear combination the following Mellin transform interpretation:
\begin{proposition}
For $\gamma_\pm\in\mathbb{R}$ we have 
\begin{align}\label{MH}
&\mathcal{M}\Bigr[\left(\gamma_+\mathcal{Z}^\lambda_++\gamma_-\mathcal{Z}^\lambda_-\right){f}\Bigr]\left(\frac{1+\Def}{2}+\frac{s}{\lambda}\right)\Biggr/\mathcal{M}[{f}]\hspace{-0.07cm}\left(\frac{1+\Def}{2}+\frac{s}{\lambda}\right)\\&=\Biggr[\hspace{-0.05cm}\frac{\gamma_++\I\gamma_{-}}{2}(1-\id)\cdot\zeta\left(\frac{(1+\Def)\lambda}{2}+s\right)+\frac{\gamma_+-\I\gamma_{-}}{2}(1-\id)\cdot\zeta\left(\frac{(1+\Def)\lambda}{2}-s\right)\hspace{-0.03cm}\Biggr]\nonumber
\end{align}
\end{proposition}
We have the implications
$(1-z)\zeta(z)=0\Rightarrow\mathcal{M}\left[\mathcal{Z}^\lambda_\pm\left(\frac{1-\lambda}{\lambda}\right){f}\right](z/\lambda)=0\;\forall\;{f}\in\mathcal{H}$ and a reformulation of section \ref{PolyEx} with help of the convolution $\ast$ , see definition \ref{MuCo}:
\begin{lemma}\label{PolyImp} The function on the r.h.s of
\begin{align*}
&\frac{4\mathcal{M}\bigr[\mathcal{Z}^\lambda_+f_1\ast\mathcal{Z}^\lambda_-f_2\bigr]\hspace{-0.15cm}\left(\frac{1+\Def}{2}+\frac{s}{\lambda}\right)}{\I\mathcal{M}[f_1\ast{f}_2]\hspace{-0.1cm}\left(\frac{1+\Def}{2}+\frac{s}{\lambda}\right)}=\bigr((1-\id)\cdot\zeta\bigr)^2\left(\frac{(1+\Def)}{2/\lambda}+s\right)-\bigr((1-\id)\cdot\zeta\bigr)^2\left(\frac{(1+\Def)}{2/\lambda}-s\right)
\end{align*}
admits only strictly complex zeros $\in\I\mathbb{R}$ if the inequality $\frac{(1+\Def)\lambda}{2}\geq1$ is satisfied
\end{lemma}
The r.h.s. of \ref{PolyImp} is by Hadamard factorization determined by its zeros while the operators $\mathcal{Z}^\lambda_\pm$ depend on their eigenvalues and eigenstates it seems plausible that there is a connection between the two relevant sets of real numbers, but without this missing link this is just a coincidence.

Notice that there would be a serious convergence problem when $\mathcal{Z}^\lambda_\pm$ gets applied to $t^\nu$, this functions are not elements of $\mathcal{H}$, the continuous spectrum of $H_\lambda\circ{\Ps^\lambda}$ does not carry over to the symmetrizations. Eigenstates $\psi_E$ of $\mathcal{Z}^\lambda_\pm$ would be a bit bizarre, we give a vague sketch how to unravel the eigenvalue equation: Either this would imply the not convincing functional equation
\begin{equation}\label{UnCo}
(1-\id)\cdot\zeta\left(\frac{(1+\Def)}{2/\lambda}+s\right)\pm(1-\id)\cdot\zeta\left(\frac{(1+\Def)}{2/\lambda}-s\right)\stackrel{!?}{=}\begin{cases}{4E\in\mathbb{R}}\\{4E/\I\in\I\mathbb{R}}\end{cases}
\end{equation}
or the function $\mathcal{M}[\psi_E]\bigr(\frac{1+\Def}{2}+\frac{s}{\lambda}\bigr)$
should be concentrated on the isolated points in $\mathbb{C}$ that are solutions of this equation \ref{UnCo}. Let us informally express the previous concentration statement in a distributional sense with help of Dirac delta functionals $\delta$ and assume
\begin{equation}\label{DeltaDi}
\mathcal{M}[\psi_E]\left(\frac{(1+\Def)}{2}+\frac{s}{\lambda}\right)=\sum{a}_z\delta(s-z)
\end{equation}
where the sum ranges over the isolated set
$${\left\{z\bigr\vert(1-\id)\zeta\left(\frac{(1+\Def)}{2/\lambda}+z\right)\pm(1-\id)\zeta\left(\frac{(1+\Def)}{2/\lambda}-z\right)=\binom{4E}{4E/\I}\right\}}$$
Under some conditions the original function can be recovered from its Mellin transform by
\begin{equation}\label{MeInv}
f(t)=\frac{1}{2\pi}\int_{\mathrm{r}-\I\infty}^{\mathrm{r}+\I\infty}\d{s}\;\mathcal{M}[f](s){t}^{-s}\end{equation}
The assumptions for this inversion are that $\mathcal{M}[f]({s})$ is analytic for $a<\Re(s)<b$ and converges absolutely in this strip, $\lim_{t\rightarrow\infty}\mathcal{M}[f](\mathrm{r}+\I{t})$ converges uniformly to zero for $a<\mathrm{r}<b$ and satisfies at any jump discontinuities the equation
$f(x)=\left[\lim_{t\rightarrow{x}^+}f(t)+\lim_{t\rightarrow{x}^-}f(t)\right]\bigr/2$.

The inverse Mellin transform \ref{MeInv} applied and the well-known formula
$\delta\bigr(g(\alpha)\bigr)=\delta\bigr(\alpha-z\bigr)/\vert{g'}(z)\vert$
at simple real zeros $z$ of $g$ immediate converts the concentration \ref{DeltaDi} into the statement that the eigenstates $\psi_E$ should be of the shape
\begin{equation}\label{DeltaNon}
\psi_E(t)=\sum_{\left\{z\bigr\vert\left((1-\id)\zeta\right)\left(\frac{(1+\Def)}{2/\lambda}+z\right)\pm\left((1-\id)\zeta\right)\left(\frac{(1+\Def)}{2/\lambda}-z\right)=\binom{4E}{4E/\I}\right\}}\frac{{a}_z}{\vert\lambda\vert}{t}^{-\left(\frac{z}{\lambda}+\frac{1+\Def}{2}\right)}
\end{equation}
Here we used in the computation of the eigenstates that the two dimensional $\delta$ distribution $\delta(z)$ can be represented by the product $\delta(x)\delta(y)$ of two one dimensional $\delta$ functions and in order to avoid integration singularities in the inverse Mellin transform it seems convenient if we blur the integration line $\mathrm{r}+\I\mathbb{R}$ to the integration strip $[\mathrm{r}-\epsilon,\mathrm{r}+\epsilon]+\I\mathbb{R}$ by averaging the integration. This blur argumentation also seems to resolve that the inverse Mellin transform is invariant under slight shifts of the integration line $\mathrm{r}+\I\mathbb{R}$. However our consideration with $\delta$ distributions is incomplete and conjectural: It is not clear that $\psi_E$ as described in \ref{DeltaNon} is a non-trivial element of $\mathcal{H}$ without specifying an infinite series $a_z\neq0$, further non-trivial convergence arguments and showing that the listed conditions for the inverse Mellin transform are not violated.

\subsubsection{\textit{Eigenstates of $\widehat{\mathcal{Z}}^\lambda_\pm$ and some products of $\zeta$ functions}}\label{EigenState2}
With the universal continuation formula \ref{PSC} we have the vanishing implications
$(1-z)\zeta(z)=0\Rightarrow\mathcal{M}[\widehat{\mathcal{Z}}^\lambda_+(\Def)f](z/\lambda)=\mathcal{M}[\widehat{\mathcal{Z}}^\lambda_-(\Def)f]\left(1+\Def-z/\lambda\right)=0\;\forall\;{f}\in\mathcal{H}$
where we substituted for the first equals $\tau_\Def^2=\id$ and for the second equals used \ref{IotaM}.

In the previous discussion the quasi eigenvalue $E=0$ and the \textit{``phase transition"} inequality
$\frac{(1+\Def)\lambda}{2}\geq1$
seems to play a special role because of the Hermite-Biehler theorem and the somehow informal interpretation that $E=0$ corresponds to zeros. In this context it may be worth to mention that we have a reformulation of \ref{PolyImp} for $\mathcal{Z}^\lambda_\pm$:
\begin{lemma}\label{Hplus}
For the sum $\widehat{\mathcal{Z}}^\lambda_+\pm\widehat{\mathcal{Z}}^\lambda_+$ and $f\in\mathcal{H}$ the implication
\begin{align}
&\frac{\mathcal{M}\left[\left(\widehat{\mathcal{Z}}^\lambda_+\pm\widehat{\mathcal{Z}}^\lambda_-\right)f\right]\hspace{-0.15cm}\left(\frac{1+\Def}{2}+\frac{s}{\lambda}\right)}{\mathcal{M}[{\tau_\Def}f]\hspace{-0.1cm}\left(\frac{1+\Def}{2}+\frac{s}{\lambda}\right)}=(1-\id)\cdot\zeta\left(\frac{(1+\Def)}{2/\lambda}\hspace{-0.05cm}+\hspace{-0.05cm}s\hspace{-0.05cm}\right)\pm(1-\id)\cdot\zeta\left(\frac{(1+\Def)}{2/\lambda}\hspace{-0.05cm}-\hspace{-0.05cm}s\hspace{-0.05cm}\right)\hspace{-0.05cm}=0\nonumber
\end{align}
$\Rightarrow{s}\in\I\mathbb{R}$ is true if $\frac{(1+\Def)\lambda}{2}\geq1$.
\end{lemma}
As we have seen in \ref{MH} the eigenstates of ${\mathcal{Z}}^\lambda_\pm$ in some sense refer to some special sums, but as we will see eigenstates of $\widehat{\mathcal{Z}}^\lambda_+$ refer to some special products, the calculations are quite analogous:

For example for eigenstates $\psi_E(\lambda)$ of $\widehat{\mathcal{Z}}^\lambda_+$ we find
\begin{align}
{\mathcal{M}[{\psi_E}]}\left(\frac{s}{\lambda}\right)&=\frac{(1-\id)\cdot\zeta(s)}{E}{\mathcal{M}[{\tau_\Def}\psi_E]}\left(\frac{s}{\lambda}\right)\\&\stackrel{!}{=}\frac{\left((1-\id)\cdot\zeta\right)(s)(1-\id)\cdot\zeta\bigr((1+\Def)\lambda-s\bigr)}{{E}^2}{\mathcal{M}[\psi_E]}\left(\frac{s}{\lambda}\right)\nonumber
\end{align}
This statement can be rewritten with $\delta$ distributions or alternative with \ref{MeInv} and the reversed formula of the convolution theorem \ref{ConvolutionTh}, namely $\mathcal{M}^{-1}[f_1\cdot{f}_2](t)=\left(\mathcal{M}^{-1}[f_1]\Conv\mathcal{M}^{-1}[{f}_2]\right)(t)$.

 If $\Re(\nu)<-1/\lambda$ then $t^\nu$ is a simultaneous eigenfunction of the commuting operators $\Ps^\lambda$ and $H_\lambda$ with the eigenvalues
$\zeta(-\lambda\nu)$
and
$1+\lambda\nu$
respectively. With this informal interpretation the $\zeta$ zeros should correspond to quasi eigenfunctions of $H_\lambda\circ\Ps^\lambda$, {\em i.e.} eigenfunctions to the eigenvalue $0$, but however for $\Re(\nu)\geq-1/\lambda$ we do not have convergence of 
$H_\lambda\circ\Ps^\lambda{t^\nu}$. If we suppose a convergent series $f(t)=\sum_{n=1}^{\infty}a_n{t}^{-n}$ we symbolically find, evaluating the operator $H_\lambda\circ\Ps^\lambda$ on each term of the expansion, the formula
\begin{equation}\label{InfoApp}
(H_\lambda\circ\Ps^\lambda{f})(t)=\sum_{n=1}^{\infty}(1-\lambda{n})\zeta(\lambda{n})a_n{t}^{-n}
\end{equation}
Another curiosity for $\widehat{\mathcal{Z}}^\lambda_\pm$ is that ${t}^{-\left(\frac{z}{\lambda}+\frac{1+\Def}{2}\right)}$ is a eigenfunction of $H_\lambda\circ\Ps^\lambda$ if convergence holds but gets mapped by ${\tau_\Def}$ to ${t}^{-\left(\frac{-z}{\lambda}+\frac{1+\Def}{2}\right)}$, hence is not respected by the compositions $\widehat{\mathcal{Z}}^\lambda_+={H}_\lambda\circ\Ps^\lambda\circ{\tau_\Def}$ and $\widehat{\mathcal{Z}}^\lambda_-={\tau_\Def}\circ{H}_\lambda\circ\Ps^\lambda$ and application of \ref{InfoApp} yields for example
\begin{align}\label{NewState}
\widehat{\mathcal{Z}}^\lambda_+\psi_E&=\sum_{\left\{z\bigr\vert(1-\id)\cdot\zeta\left(\frac{(1+\Def)}{2/\lambda}+z\right)\cdot(1-\id)\cdot\zeta\left(\frac{(1+\Def)}{2/\lambda}-z\right)={E^2}\right\}}\hspace{-2cm}(1-\id)\cdot\zeta\left(\frac{1+\Def}{2/\lambda}+z\right)\frac{{a}_z}{\vert\lambda\vert}{t}^{-\left(\frac{-z}{\lambda}+\frac{1+\Def}{2}\right)}\\&=\sum_{\left\{z\bigr\vert(1-\id)\cdot\zeta\left(\frac{(1+\Def)}{2/\lambda}+z\right)\cdot(1-\id)\cdot\zeta\left(\frac{(1+\Def)}{2/\lambda}-z\right)={E^2}\right\}}E\frac{{a}_z}{\vert\lambda\vert}{t}^{-\left(\frac{z}{\lambda}+\frac{1+\Def}{2}\right)}
\end{align}
This implies by comparison of the two expressions for $\widehat{\mathcal{Z}}^\lambda_+\psi_E$ the equation $E\frac{{a}_{-z}}{a_z}=(1-\id)\cdot\zeta\left(\frac{1+\Def}{2/\lambda}+z\right)$, hence for the coefficients $a_z$ the relation
$\pm\frac{{a}_{-z}}{a_z}=\sqrt{\frac{(1-\id)\cdot\zeta\left(\frac{1+\Def}{2/\lambda}+z\right)}{(1-\id)\cdot\zeta\left(\frac{1+\Def}{2/\lambda}-z\right)}}$.

As will be alluded in section \ref{ConLem} analogous statements also hold for the eigenstates of the four standard symmetrizations of an arbitrary convolution operator obtained by lemma \ref{SACon}.

Compared with the general $E\rightarrow-E$ eigenvalue symmetry construction \ref{InvoCon} with help of $\tau_\Def$ we have a different $E\rightarrow-E$ observation for $\widehat{\mathcal{Z}}^\lambda_\pm$ because $H_\alpha$ commutes with $\Ps^\lambda$ and \ref{ComRel1}:
\begin{proposition}\label{ImpiLaPla}
The following implications are true
\begin{align*}
\widehat{\mathcal{Z}}^\lambda_\pm\psi_E=E\psi_E&\Rightarrow\widehat{\mathcal{Z}}^\lambda_\pm\left(H_{\frac{2}{1+\Def}}\psi_E\right)=-E\left(H_{\frac{2}{1+\Def}}\psi_E\right)\nonumber\Rightarrow\widehat{\mathcal{Z}}^\lambda_\pm\left(\Delta_{\frac{2}{1+\Def}}\psi_E\right)=E\left(\Delta_{\frac{2}{1+\Def}}\psi_E\right)
\end{align*}
\end{proposition}
It is therefore immediate to consider the quotient of the eigenspace of the eigenvalue $E$ by the equivalence relation
$$\psi_E\sim\psi_E'\Leftrightarrow\exists{n}:\psi_E=\Delta_{\frac{2}{1+\Def}}^n\psi_E'\vee\psi_E'=\Delta_{\frac{2}{1+\Def}}^n\psi_E$$

With \ref{InfoApp} $H_{\frac{2}{1+\Def}}$ acts on the coefficients $a_z$ by multiplication with
$-\frac{2z}{\lambda(1+\Def)}$
while the second order differential operator $\Delta_{\frac{2}{1+\Def}}$ acts by multiplication with
$\bigr(\bigr(\frac{2z}{\lambda(1+\Def)}\bigr)^2-1\bigr)$.

The case $(1+\Def)\lambda=1$ is special: Here we have that the eigenstate $\psi_0$ has to be of the shape $\sum_{z\vert\zeta(1/2+z)=0}{t}^{-\left(\frac{z}{\lambda}+\frac{1+\Def}{2}\right)}$  and Weil showed \cite{We} that RH is equivalent to the remarkable positivity criteria that for Mellin invertible test functions $f$
\begin{equation}\label{Weil}
\sum\mathcal{M}\left[{\tau_0}\overline{f}\Conv{f}\right](z)\geq0
\end{equation}
where the sum runs over zeros $z$ of ${\Xi}$ and $\Conv$ denotes the multiplicative convolution, see definition \ref{MuCo}. We refer the reader for proofs and a detailed discussion of  Weil's positivity criteria \ref{Weil} and  his so-called explicit formula of number theory to \cite{Bur}.
\section{\textit{The analogy with convolution operators}}\label{ConLem}
The algebra of multiplication operators \ref{mult} became singular because mild multiplied singularities at $0$ or $\infty$ get absorbed by the nice behaved test functions in $\mathcal{H}$. For the multiplication operators $\Mu_V$ it seems not so easy to give a general Mellin transform interpretation we only have by integration by parts the symmetric formula
\begin{align}\label{MHMu}
s\mathcal{M}\left[{V}\cdot{f}\right]\left(s\right)&=s\mathcal{M}\left[\Mu_{V}{f}\right]\left(s\right)=s\mathcal{M}\left[\Mu_{f}{V}\right]\left(s\right)=-\mathcal{M}\left[\Mu_{t\partial_t{V}}f\right]\left(s\right)-\mathcal{M}\left[\Mu_{t\partial_t{f}}V\right]\left(s\right)\end{align}

From this interpretation point of view there is a more well-behaved algebra acting on $\mathcal{H}$ by integration:
Let as usual the multiplicative convolution $\Conv$ be defined by 
\begin{equation}\label{MuCo}
(f_1\Conv{f}_2)(t):=\int_0^\infty\frac{\d{y}}{y}f_1\left(\frac{t}{y}\right)f_2(y)
\end{equation}
The bilinear operation $\ast$ is associative, commutative and we have the famous convolution theorem
\begin{equation}\label{ConvolutionTh}
\mathcal{M}[f_1\Conv{f}_2](s)=\mathcal{M}[f_1](s)\cdot\mathcal{M}[{f}_2](s)
\end{equation}
With the usual additive convolution
$(f_1\tilde{\Conv}{f}_2)(t):=\int_{-\infty}^\infty{\d{y}}f_1\left({t}-{y}\right)f_2(y)$
we recover $\Conv$ by a so-called equivalence formula that uses invertibility to make a trivial associativity connection by the identity ${(f}\tilde{\Conv}{g})\left({t}\right)=\bigr((f\circ\ln)\Conv({g}\circ\ln)\bigr)\left(\exp(t)\right)$.

For the norm $\vert\vert\cdot\vert\vert_p$ defined by $\vert\vert{f}\vert\vert_p:=\bigr(\int_{-\infty}^\infty\d{x}\vert{f}(x)\vert^p\bigr)^{1/p}$ the Young inequality
$\vert\vert{{f}\tilde{\Conv}{g}}\vert\vert_r\leq\vert\vert{f}\vert\vert_p\vert\vert{g}\vert\vert_q$ is valid for $\frac{1}{p}+\frac{1}{q}=1+\frac{1}{r}$ where $1\leq{p},{q},{r}\leq\infty$.

 It is also well-known that functions of the shape ${f_1}\ast{f_2}$ are dense in $\mathcal{H}$.

\begin{definition}\label{MuCon} Let $\Co_V:\mathcal{H}\rightarrow\mathcal{H}$ with fixed $V\in\mathcal{H}$ act by
$f(t)\rightarrow\left(\Co_V{f}\right)(t):=({V}\Conv{f})(t)$.
\end{definition}
We have for example
$\Co_{\delta(1-t)}=\id$
where $\delta$ denotes as usual the Dirac delta function, but notice $\delta\notin\mathcal{H}$ and $\delta$ can only be approximated elements of $\mathcal{H}$. Another interesting example of an identity for the convolution operator $\Co_{e^{-\rho\ln^2(t)}t^s}$ is the interpretation
\begin{equation}\label{e2Conv}
\left(\Co_{e^{-\rho\ln^2(x)}x^s}f\right)(t)={e^{-\rho\ln^2(t)}t^s}\mathcal{M}\left[f{e^{-\rho\ln^2(\cdot)}}\right]\Bigr(-s+2\rho\ln(t)\Bigr)
\end{equation}
The ring $(\mathcal{H},+,\Conv)$ is clearly not noetherian, for instance the spaces $\Mu_{e^{-\rho\ln^2(t)}}\mathcal{H}$ where $0\leq\rho$ are $\Conv$-ideals and the inclusion
$\Mu_{e^{-\rho'\ln^2(t)}}\mathcal{H}\subset\Mu_{e^{-\rho\ln^2(t)}}\mathcal{H}$
holds if $\rho<\rho'$.
\begin{proposition}\label{CoAd}
With respect to the inner product \ref{SP} the adjoint of $\Co_V$ is given by
$$\Co_V^\ast=\mathrm{Con}_{\tau_\Def}\Co_{\overline{V}}={\tau_\Def}\circ\Co_{\overline{V}}\circ\tau_\Def=\Co_{{\tau_\Def}\overline{V}}$$
\end{proposition}

By the convolution theorem \ref{ConvolutionTh} we have for $f,V\in\mathcal{H}$ the Mellin transform interpretation
$\mathcal{M}\left[\Co_V{f}\right]\left(s\right)=\mathcal{M}\left[\Co_f{V}\right]\left(s\right)=\mathcal{M}\left[V\right]\left(s\right)\cdot\mathcal{M}\left[{f}\right]\left(s\right)$,
especially 
$\mathcal{M}\left[\Co_{\Delta_4\Psi}{f}\right]\left(s\right)=4\xi(s)\mathcal{M}\left[{f}\right]\left(s\right)$
hence \ref{MuCon} allows to encode $\zeta$ alternative to \ref{Ad-Hoc} as a convolution Hamiltonian:
\begin{lemma}\label{SACon} For every $V\in\mathcal{H}$ we have $\Co_{(V+{\tau_\Def}\overline{V})/2}\in\mathfrak{s}$ and $\Co_{\I(V-{\tau_\Def}\overline{V})/2}\in\mathfrak{s}$. If $V(t)=\overline{V(t)}$ we have also two self-adjoint operators given by the compositions
${\tau_\Def}\circ\Co_{V}\in\mathfrak{s}$ and $\Co_{V}\circ{\tau_\Def}\in\mathfrak{s}$.\end{lemma}
\begin{proof}[Proof] We use \ref{ScaCon}, the fact that $\Conv$ is associative and the compatibility
\begin{equation}\label{CoCo1}
{\tau_\Def^\mu}(f_1\Conv{f}_2)={\vert\mu\vert}({\tau_\Def^\mu}{f}_1\Conv{\tau_\Def^\mu}{f}_2)
\end{equation}
This shows in particular that the multiplication operators $\Mu_{t^s}$ are morphisms of $\Conv$, {\em i.e.} $\Mu_{t^s}\circ\Conv=\Conv\circ(\Mu_{t^s}\otimes\Mu_{t^s})$ and more general
$\tau^\mu_\Def/\vert\mu\vert$ are morphisms of $\Conv$. We have the commutation relation
\begin{equation}\label{CoHam}
{\tau_\Def^\mu}\circ\Co_V=\vert\mu\vert\Co_{{\tau_\Def^\mu}{V}}\circ{\tau_\Def^\mu}
\end{equation}
The first two Hamiltonians are now obtained by the two standard symmetrization procedures $\mathrm{p}$ and $\partial$ and the third and fourth by the standard procedures described in \ref{proof2}.\end{proof}
 Notice that there is literally in some sense an analogy of the following four symmetrizations in the lemma \ref{SACon} with the four more discrete formulas \ref{Ad-Hoc} defined by summation.


With ${\tau_\Def}$, $\Di_t$ and $\langle\cdot,\cdot\rangle$ we can rewrite the convolution $\Conv$ in the shape
\begin{equation}\label{ScaCon}
({\tau_\Def}\overline{f_1}\Conv{f}_2)(t)=\langle{\tau_\Def}\circ\Di_{t}\circ{\tau_\Def}{f_1},f_2\rangle_\Def=\langle{f_1},\Di_{t}f_2\rangle_\Def
\end{equation}
With the reformulation \ref{ScaCon} and the compatibilities \ref{CoCo1}, \ref{CoCo3} and $\tau_\Def\circ\Di_t=\Mu_{t^{-(1+\Def)}}\Di_t\circ\tau_\Def$ we get a reincarnation of theorem \ref{Ad-Hoc} in a simple, basic form concerning convolution operators:
\begin{lemma}\label{CoAd-Hoc} We have the compatibility
$$H_\lambda\circ\Ps^\lambda(f_1\Conv{f}_2)=(H_\lambda\circ\Ps^\lambda{f_1})\Conv{f}_2={f_1}\Conv({H}_\lambda\circ\Ps^\lambda{f}_2)$$
We have $\mathcal{Z}^\lambda_\pm\circ\Co_V\in\mathfrak{s}$ if $\overline{{\tau_\Def}(V)}=V$ and $\widehat{\mathcal{Z}}^\lambda_\pm\circ\Co_V\in\mathfrak{s}$ if $V$ is real.
\end{lemma}
\begin{proof}[Proof]
By \ref{CoCo1} and the standard identities
\begin{equation}\label{CoCo2}
H_\alpha(f_1\Conv{f}_2)=H_\alpha{f}_1\Conv{f}_2=(H_\alpha{f}_1\Conv{f}_2+{f}_1\Conv{H}_\alpha{f}_2)/2
\end{equation}
\begin{equation}\label{CoCo3}
{\Di_\beta}(f_1\Conv{f}_2)={\Di_\beta}{f}_1\Conv{f}_2=(\Di_\beta{f}_1\Conv{f}_2+{f}_1\Conv\Di_\beta{f}_2)/2
\end{equation}
Hence $\mathcal{Z}^\lambda_\pm\circ\Co_V$ is self-adjoint if $\overline{{\tau_\Def}(V)}=V$ because we have the identity
\begin{equation}
\mathcal{Z}^\lambda_\pm\circ\Co_V=\Co_V\circ\mathcal{Z}^\lambda_\pm
\end{equation}
By \ref{CoHam} if
$\mathcal{Z}^\lambda_\pm\psi_E=E\psi_E$
we have the eigenvalue equation
$\mathcal{Z}^\lambda_\pm(\Co_{V}\psi_{E})=E(\Co_V\psi_{E})$, the eigenstates of $\mathcal{Z}^\lambda_\pm$ are stable under convolution and the convolution of eigenstates with different eigenvalues vanishes. The eigenvalues of $\widehat{\mathcal{Z}}_\pm^\lambda$ are not stable under convolution in general: 
\begin{equation}\label{HCONV}
\widehat{\mathcal{Z}}_\pm^\lambda\circ\Co_V=\Co_{{\tau_\Def}{V}}\circ\widehat{\mathcal{Z}}_\pm^\lambda
\end{equation}
holds, hence if
$\widehat{\mathcal{Z}}_\pm^\lambda\psi_E=E\psi_E$
we find
$\widehat{\mathcal{Z}}_\pm^\lambda(\Co_{V}\psi_{E})=E(\Co_{{\tau_\Def}{V}}\psi_{E})$
and this implies for potential $V^\pm\in\mathcal{H}^\pm_\Def$ in that
$\widehat{\mathcal{Z}}_\pm^\lambda(\Co_{V^\pm}\psi_{E})=\pm{E}(\Co_{V^\pm}\psi_{E})$. This shows that $\widehat{\mathcal{Z}}^\lambda_\pm\circ\Co_V\in\mathfrak{s}$ if $V$ is real.\end{proof}
The equation $\mathcal{M}\left[\Mu_{\ln(t)}f\right]\left(s/\lambda\right)=\lambda\partial_s\mathcal{M}\left[f\right]\left(s/\lambda\right)$ holds, hence we have the identity
$$\mathcal{M}\left[\underbrace{[\cdots[}_{n}\Ps^\lambda,\Mu_{\ln(t)}]_-\cdots,\Mu_{\ln(t)}]_-{f}\right]\left(s/\lambda\right)=\lambda^n\left(\partial^n_s\zeta(s)\right)\mathcal{M}\left[{f}\right]\left(s/\lambda\right)$$
Also a similar compatibility is satisfied for the commutators with convolution operators: The fact that $\Mu_{\ln(t)}$ is a derivation of $\Conv$, {\em i.e.}
$0=\Conv\circ(\id\otimes\Mu_{\ln(t)}+\Mu_{\ln(t)}\otimes\id)-\Mu_{\ln(t)}\circ\Conv:\mathcal{H}^{2\times}\rightarrow\mathcal{H}$, is a gadget in the proof of the explicit formulas of number theory contained in \cite{RM}. This derivation property and the commutativity of $\Conv$ allows to define a Jacobi structure:
\begin{lemma}\label{DerLie} Let $(\mathcal{A}^\bullet,\star)$ be a graded super-algebra and $\partial:\mathcal{A}^\bullet\rightarrow\mathcal{A}^{\bullet+2m}$ with $m\in\mathbb{Z}$ be an even super-derivation {\em i.e.} $\star$ is associative, $\vert{a\star{b}}\vert=\vert{a}\vert+\vert{b}\vert$ and $\partial(a\star{b})=(\partial{a})\star{b}+a\star(\partial{b})$. Suppose $\mathcal{A}$ is super-commutative {\em i.e.} $a\star{b}=(-1)^{\vert{a}\vert\vert{b}\vert}b\star{a}$ or that $\partial$ is a differential {\em i.e.} $\partial^2=0$: By $[x,y]_\partial:=(\partial{x})\star{y}-(-1)^{\vert{x}\vert\vert{y}\vert}(\partial{y})\star{x}$ we have a even super Lie-bracket $[\cdot,\cdot]_\partial:\mathcal{A}^{\otimes2}\rightarrow\mathcal{A}$ of degree $2m$, {\em i.e.} $(-1)^{\vert{x}\vert+\vert{z}\vert}[x,[y,z]_\partial]_\partial+(-1)^{\vert{y}\vert+\vert{x}\vert}[y,[z,x]_\partial]_\partial+(-1)^{\vert{z}\vert+\vert{y}\vert}[z,[x,y]_\partial]_\partial=0$.
\end{lemma}
For the usual derivative we have the interpretation of $[x,y]_\partial=y^2\partial(x/y)$.

\section{\textit{Substitution Hamiltonians and uniqueness of ${\tau_\Def}$}}\label{SubSti}
\begin{definition}\label{UpG}
Suppose $g:[0,\infty]\rightarrow[0,\infty]$ is invertible, differentiable and respects rapid decay at $0$ and $\infty$ in the sense that if $f$ decays rapidly at this two points then also $f\circ{g}^{(-1)}$ have this property and let $V$ be as specified in \ref{MuLti}.
\begin{equation}\label{Adju1}
\left(\mathcal{S}_{V,g}f\right)(t):=\bigr(\left(\Mu_{{V}}f\right)\circ{g}\Bigr)(t)={V\left(g(t)\right)}f\left(g(t)\right)
\end{equation}
\end{definition}
Definition \ref{defiota} is an example of a substitution operator, {\em i.e.} we have ${\tau_\Def^\mu}=\mathcal{S}_{t^{1+\Def},t^\mu}$ and also the dilation can be written as $\Di_\beta=\mathcal{S}_{1,\beta{t}}$. Clearly \ref{UpG} also unifies with \ref{MuLti} for $g=\id$. For the composition of two substitution operators we have the semi-direct composition rule
\begin{equation}
\mathcal{S}_{V_1,g_1}\circ\mathcal{S}_{V_2,g_2}=\mathcal{S}_{(V_1\circ{g_2^{-1}})\cdot{V_2},g_2\circ{g}_1}
\end{equation}

The following proposition at least allows us by the mentioned procedures $\mathcal{S}_{V,g}\rightarrow\mathrm{p}\mathcal{S}_{V,g}\in\mathfrak{s}$ and $\mathcal{S}_{V,g}\rightarrow\partial\mathcal{S}_{V,g}\in\mathfrak{s}$ {\em etc.}
to produce some self-adjoint substitution operators or use the conjugation
$A\rightarrow\mathrm{Con}_{\mathcal{S}_{V,g}}^\ast(A)={\mathcal{S}_{V,g}}\circ{A}\circ\mathcal{S}_{V,g}^\ast$
to transform $A\in\mathfrak{s}$.
\begin{lemma}\label{SuAd}
The adjoint of the operator $\mathcal{S}_{V,g}$ with respect to \ref{SP} is given by
\begin{equation}\label{Adju2}
\left(\mathcal{S}^\ast_{V,g}f\right)(t)=\pm\Biggr(\mathcal{S}_{\frac{t^\Def\cdot(\overline{V}\circ{g})}{g^\Def\cdot\partial_t{g}},g^{-1}}f\Biggr)(t)=\pm\partial_t{g}^{-1}(t)\overline{{V}(t)}\left(\frac{g^{-1}(t)}{t}\right)^\Def{f}\left(g^{-1}(t)\right)
\end{equation}
with the $+$ sign in the case $g(0)=0$ and the $-$ sign in the case $g(0)=\infty$.
\end{lemma}
We have the following corollary of \ref{Adju2} that not surprisingly just identifies two times the same substitution operator on both sides of $\langle\cdot,\cdot\rangle_\Def$ with a multiplication operator:
\begin{korollar}\label{SubAdjo}
Let $g(0)=0$. We have for $f_1,f_2\in\mathcal{H}$ the adjunction equation
$$\left\langle\mathcal{S}_{V,g}{f_1},\mathcal{S}_{V,g}{f_2}\right\rangle_\Def=\pm\left\langle{f_1},\vert{V(t)}\vert^2\partial_t{g}^{-1}(t)\left(\frac{g^{-1}(t)}{t}\right)^\Def{f_2}\right\rangle_\Def$$
Hence we have $\big\langle\mathcal{S}_{V,g}{f_1},\mathcal{S}_{V,g}{f_2}\big\rangle_\Def=\left\langle{f_1},\Mu_{h}{f_2}\right\rangle_\Def$ with some real positive $h\in\mathcal{H}$ if we set
$V(t)=e^{\I{\varphi(t)}}\sqrt{\left(\frac{t}{g^{-1}(t)}\right)^\Def\frac{h{(t)}}{\partial_t{g}^{-1}(t)}}$
with some real argument $\varphi:\mathbb{R}^+\rightarrow\mathbb{R}$.
\end{korollar}
 Clearly we also have compatibility of the Lie-algebra structure and the \textbf{Chevalley-Eilenberg complex} with  
adjoint conjugation by $\mathrm{Con}^\ast_{\mathcal{S}_{V,g}}$ as specified in \ref{SubAdjo}.

If we scale ${\tau_\Def^\mu}$ with $\sqrt{\vert\mu\vert}$ we could omit the proportionality factor $1/\vert\mu\vert$ in the adjunction $\bigr({\tau_\Def^\mu}\bigr)^\ast={{\tau_\Def^{1/\mu}}}/{\vert\mu\vert}$, {\em i.e.} $(\sqrt{\vert\mu\vert}{\tau_\Def^\mu})^\ast=\sqrt{\vert1/\mu\vert}{\tau_\Def}^{1/\mu}$
and the previous observation is a solution of a question that we only answer halfway in the following proposition \ref{PartiAnsw}:

We ask for the symmetry that \ref{Adju2} is $1/\mu$ proportional to the operator defined in \ref{Adju1} but with $g$ replaced by its inverse $g^{-1}$, {\em i.e.} we search for real pairs $(V,g)$ with 
$\sqrt{\mu}\mathcal{S}^\ast_{V,g}=\sqrt{1/{\mu}}\mathcal{S}_{V,g^{-1}}$.
\begin{proposition}\label{PartiAnsw}
If $g$ satisfies 
\begin{equation}\label{Consti1}
\mu^{\mu^2-1}\partial_t{g}\left(g(t)\right)=\left(\partial_t{g}^{-1}(t)\right)^{-\mu^2}=\left({\partial_t{g}\left(g^{-1}(t)\right)}\right)^{\mu^2}
\end{equation}
then we have $\sqrt{\mu}\mathcal{S}^\ast_{V,g}=\sqrt{1/{\mu}}\mathcal{S}_{V,g^{-1}}$ if $V(t)$ is real proportional to
\begin{equation}\label{FunSol}
\mu^{\frac{2}{1-\mu^2}}t^\Def\left(\frac{1}{\partial_t{g}(t)\partial_t{g}\left(g(t)\right)}\right)^{\frac{1}{1-\mu^2}}
\end{equation}
\end{proposition}
\begin{proof}[Proof] In concrete formulas we ask for
\begin{align}\label{Adju3}
&\frac{{t}^\Def}{V(t)}=\mu\frac{\left(g^{-1}(t)\right)^\Def}{{V}\left(g^{-1}(t)\right)}\partial_t{g}^{-1}(t)\quad\forall{t}\in[0,\infty]
\end{align}
For some $V$ a $g$ solving \ref{Adju3} is unique by the Picard-Lindel\"of theorem. Iteration of \ref{Adju3} and the Ansatz $V(t)=t^\Def{v}(t)$ decouples the parameter $\Def$ and is now more symmetric in $g$ and $g^{-1}$:
\begin{equation}\label{NewFun}
\mu^2\frac{\partial_t{g}^{-1}(t)}{v\left(g^{-1}(t)\right)}=\mu^2\partial_t\ln\left({v\left(g^{-1}(t)\right)}\right)=\partial_t\ln\bigr({v\left(g(t)\right)}\bigr)=\frac{\partial_t{g}(t)}{v\left(g(t)\right)}
\end{equation}
By integration and $g(0)=g^{-1}(0)$ we yield
${v\left(g(t)\right)}={v^{1-\mu^2}\bigr(g(0)\bigr)}{v^{\mu^2}\left(g^{-1}(t)\right)}$.
Self-consistency with \ref{NewFun} and the more exquisite equation \ref{Adju3} implies if $\mu^2\neq1$ the constraint
\ref{Consti1}.
If $g$ satisfies \ref{Consti1} we have
$\sqrt{\mu}\mathcal{S}^\ast_{V,g}=\sqrt{1/{\mu}}\mathcal{S}_{V,g^{-1}}$
if $V(t)$
is real proportional to
\ref{FunSol}. \end{proof}
\begin{lemma}\label{SubInvo}
Let $g$ diverge at $0$ like some power $t^{-\rho}$ with $\rho\in\mathbb{R}^+$. If the graph of $g$ is reflection symmetric with respect to the diagonal we have an involutive operator $\mathcal{S}_{V,g}$ where
$V=\sqrt{\frac{v}{{v}\circ{g}}}$
with some $v:\mathbb{R}^+\rightarrow\mathbb{R}^+$ with no zeros and so that $V$ is as specified in \ref{MuLti}. 
\end{lemma}\begin{proof}[Proof] If we suppose $\mathcal{S}_{V,g}$ is an involution we obtain $g\circ{g}=\id$ and $V(t)\cdot{V}(g(t))=1$.
Notice $g\circ{g}=\id$ {\em i.e.} $g=g^{-1}$ implies if $g(0)=0$ readily $g=\id$ but else is locally not very restricting, for instance it is well-known that the inverse of $g$ is just the graph obtained by reflection at the dashed diagonal in the picture below:
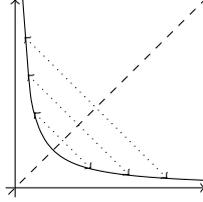
\begin{center}
\vspace{-0.5cm}\begin{tikzpicture}
  \draw[->] (-0.125,0) -- (2.5,0) node[right] {};
  \draw[dotted,<->] (0.125,2) -- (2,0.125);
  \draw[dotted,<->] (0.25/1.5,1.5) -- (1.5,0.25/1.5);
  \draw[dotted,<->] (0.25,1) -- (1,0.25);
  \draw[->] (0,-0.125) -- (0,2.5) node[above] {};
  \draw[dashed, scale=0.5,domain=-0.25/sqrt(2):5,smooth,variable=\x,black] plot ({\x},{\x});
  \draw[scale=0.5,domain=0.2:5,smooth,variable=\y,black] plot ({\y},{1/(\y)});
\end{tikzpicture}\vspace{-0.8cm}\end{center}\end{proof}
Hence there are plenty involutive substitution operators, but self-adjoint ones are quite rare and the two symmetries of the hero $\tau_\Def$ of the previous section \ref{Std1} in some sense are unique:
\begin{proposition}
If we suppose a substitution operator $\mathcal{S}_{V,g}$ as specified in \ref{UpG} satisfies the two identities $\mathcal{S}_{V,g}\circ\mathcal{S}_{V,g}=\id$ and $\sqrt{\mu}\mathcal{S}^\ast_{V,g}=\sqrt{1/{\mu}}\mathcal{S}_{V,g^{-1}}$ then $\mathcal{S}_{V,g}=\id$ or $\mathcal{S}_{V,g}={\tau_\Def}$.
\end{proposition}
\begin{proof}[Proof] We find $\mathcal{S}_{V,g}=\id\vee{\tau_\Def}$ by substituting  $g\circ{g}=\id$ and $V(t)\cdot{V}(g(t))=1$ in \ref{FunSol}.\end{proof}
If $B^\ast=B^{-1}$ the adjoint conjugation $\mathrm{Con}^\ast_{B}(A)=B\circ{A}\circ{B}^\ast$ obviously coincides with the usual conjugation $\mathrm{Con}_{B}(A):=B\circ{A}\circ{B}^{-1}$ but $\mathrm{Con}^\ast_{B}$ restricts to $\mathfrak{s}$ $\forall{B}\in\mathfrak{g}$ because $\mathrm{Con}^\ast_{B}$ commutes with $^\ast$ while $\mathrm{Con}_{B}$ is only defined if $B$ is invertible in the sense $B\in\mathfrak{i}$
$\mathfrak{i}:=\{A\in\mathfrak{g}\vert\exists!{A}^{-1}\in\mathfrak{g}:A\circ{A}^{-1}=A^{-1}\circ{A}=\id\}$.
Let us mention that we have a map $\mathcal{Q}:\mathfrak{i}\rightarrow\mathfrak{i}\cap\mathfrak{s}(\Def)$ by $\mathcal{Q}(A):={A}\circ{A}^{\ast(\Def)}$ and because $(A^{-1})^\ast=(A^\ast)^{-1}$ we have $\mathcal{Q}\bigr((A^\ast)^{-1}\bigr)=\left(\mathcal{Q}(A)\right)^{-1}$. The map $\mathrm{Con}_{B}$ restricts to $\mathfrak{s}$ if $B=B^{-1}\in\mathfrak{s}$ and clearly we have $\mathrm{Con}_{A}\circ\mathrm{Con}_{B}=\mathrm{Con}_{A\circ{B}}$ as well as $\mathrm{Con}^\ast_{A}\circ\mathrm{Con}^\ast_{B}=\mathrm{Con}^\ast_{A\circ{B}}$.

\begin{proposition}
For the set of smooth involutions of $\mathbb{R}^+$ the quotient by the equivalence relation $\sim$ defined by $g_1\sim{g}_2\Leftrightarrow\exists{f}:\mathrm{Con}_f{g_1}=g_2$, where $f$ is a smooth function on $\mathbb{R}^+$, just consists of the maps $\id$ and $\frac{1}{\id}$.\end{proposition}
\begin{proof}[Proof] Let $g$ be an involution with $g(0)=\infty$ and $h:\mathbb{R}^+\rightarrow\mathbb{R}^+$ with $h(0)=0$ and $h'(t)>0$ set $f(t):=h(t)/h\left(g(t)\right)$. The function $f:\mathbb{R}^+\rightarrow\mathbb{R}^+$ is bijective because $f(0)=0$, $f(\infty)=\infty$ and $f'(t)>0$. We also clearly have $1/f(t)=f\left(g(t)\right)$ or equivalent $\mathrm{Con}_f{g}=\frac{1}{\id}$.\end{proof}

\subsection{\textit{Basic compatibilities with uncertainty}}
\begin{definition}\label{Variance} The variance $\sigma_f(A)$ of $A\in\mathfrak{s}$ is defined by
$\sigma_f(A):=\sqrt{\Big\langle{f},\bigr(A-\big\langle{f,A{f}}\big\rangle_\Def\bigr)^2f\Big\rangle_\Def}$
where $f\in\mathcal{H}$ is a normalized state, i.e $f$ is on the unit sphere defined by the condition $\langle{f},{f}\rangle_\Def=1$.\end{definition}
Only the quotient space $\mathfrak{g}/\sim_\Def$ with respect to the equivalence relation
${A}\sim_\Def{B}\Leftrightarrow\exists{\mathcal{S}_{V,g}}\in\mathfrak{g}:\mathcal{S}_{V,g}\circ\mathcal{S}_{V,g}^{\ast(\Def)}=\id\;\wedge\;{A}=\mathrm{Con}^{\ast(\Def)}_{\mathcal{S}_{V,g}}B$
is relevant for $\langle\cdot,\cdot\rangle_\Def$ restricted to the unit sphere. 
\begin{korollar}\label{ConInvo}
If $V$ satisfies $\vert{V(t)}\vert^2=\left({t}/{g^{-1}(t)}\right)^\Def\bigr/\partial_t{g^{-1}}(t)$ then $\mathcal{S}_{V,g}^\ast\circ\mathcal{S}_{V,g}=\id$. This operators form a subgroup of $\mathfrak{g}$ and $\mathrm{Con}^\ast_{\mathcal{S}_{V,g}}$ is compatible with the natural Lie-algebra structure.
\end{korollar}
It is well-known that the Cauchy-Schwarz inequality
$\langle{f,f}\rangle_\Def\langle{g,g}\rangle_\Def\geq\vert\langle{f,g}\rangle_\Def\vert^2$
implies
\begin{equation}\label{SchroeUncer}
{\sigma^2_f(A_1)\sigma^2_f(A_2)\geq\Big\vert\Bigr\langle{f},\frac{1}{2\I}[A_1,A_2]_{-}f\Bigr\rangle_\Def\Big\vert^2+\Big\vert{\Bigr\langle{f},\frac{1}{2}[A_1,A_2]_{+}f\Bigr\rangle_\Def-\langle{f},A_1{f}\rangle_\Def\langle{f},A_2{f}\rangle_\Def}\Big\vert^2}
\end{equation}
for self-adjoint operators $A_i\in\mathfrak{s}$ called \textit{\textbf{Schr\"odinger uncertainty relation}} and this slightly stronger inequality clearly implies the usual \textit{Heisenberg uncertainty principle}.
\begin{proposition} Let $f=f^++f^-$ with $f^\pm\in\mathcal{H}^\pm_\Def$. For $A\in\mathfrak{s}(\Def)$ we have
$$\Bigr\langle{f,\frac{1}{2}[{A},{\tau_\Def}]_{\pm}f}\Bigr\rangle_\Def=2\Re\left(\Bigr\langle{f^{-}},\frac{1}{2}[{A},{\tau_\Def}]_{\pm}{f^{+}}\Bigr\rangle_\Def\right)+\langle{f^\pm,A{f}^\pm}\rangle_\Def-\langle{f^-,A{f}^-}\rangle_\Def$$
\end{proposition}
\begin{proof}[Proof] The variance of ${\tau_\Def}$ in a normalized eigenstate ${f^{\pm}}/{\sqrt{\langle{f^{\pm},f^{\pm}}\rangle_\Def}}$ vanishes, hence the l.h.s. of \ref{SchroeUncer}. For self-adjoint $B$ the expectation value
$\langle{g,{B}g}\rangle_\Def=\overline{\langle{{B}g,g}\rangle_\Def}=\overline{\langle{g,{B}g}\rangle_\Def}$
is real, hence the square greater than zero if the expectation value is $\neq0$. The resulting equations are invariant under scaling, we can drop the unit sphere normalization and yield
$\left\langle{f^{\pm}},\frac{\I}{2}[{A},{\tau_\Def}]_{-}f^{\pm}\right\rangle_\Def=0=\left\langle{f^{\pm}},\left(\frac{1}{2}[{A},{\tau_\Def}]_{+}\mp{A}\right){f^{\pm}}\right\rangle_\Def$.
The proposition is now obtained by multiplying out this in $f\in\mathcal{H}$ quadratic expressions and substitution of 
$\overline{\langle{f,Ag}\rangle_\Def}=\langle{Ag,f}\rangle_\Def=\langle{g,A^\ast{f}}\rangle_\Def$.
\end{proof}
 It is standard that the direct sum $\oplus_{i\in{I}}\mathcal{H}_i$ of an indexed family  $(\mathcal{H}_i,\langle\cdot,\cdot\rangle_i)_{i\in{I}}$ of Hilbert spaces becomes a Hilbert space with inner product $\langle\cdot,\cdot\rangle$ on $\oplus_{i\in{I}}\mathcal{H}_i$ specified by
$\left\langle{f},{g}\right\rangle:=\sum_{i\in{I}}\langle{f_i},{g_i}\rangle_i$
and this sum converges because it is by definition finite, the same holds for pre-Hilbert spaces.
\begin{proposition}
Let $f\in(\mathcal{H},\langle\cdot,\cdot\rangle_\Def)_{\Def\in(-\infty,\infty)}$ be normalized with support contained in $[0,1]$, $[1,\infty]$ or the union of two intervals $[a_{1,2},b_{1,2}]$ with $0<{a_1\leq{b_1}}\leq1\leq{a_2\leq{b_2}}<\infty$ and $a_1a_2\geq1$ or $1\geq{b_1b_2}$. For every self-adjoint $A:\oplus_{\Def\in(-\infty,\infty)}\mathcal{H}\rightarrow\oplus_{\Def\in(-\infty,\infty)}\mathcal{H}$ we have the inequality
\begin{equation}\label{SchroeUncer2}
\sigma^2_f(A)\geq\Bigr\vert\Bigr\langle{f},\frac{1}{2\I}[A,{\tau_\Def}]_{-}f\Bigr\rangle\Bigr\vert^2+\Bigr\vert{\Bigr\langle{f},\frac{1}{2}[A,{\tau_\Def}]_{+}f\Bigr\rangle}\Bigr\vert^2\hspace{-0.1cm}\end{equation}
In particular \ref{SchroeUncer} for $A_1=\widehat{\mathcal{Z}}^{\lambda_\Def}_+$ and $A_2={\tau_\Def}$ if $\mathrm{supp}(f_\Def)\subset[1,\infty]$ or for $A_1=\widehat{\mathcal{Z}}^{\lambda_\Def}_-$ with $\lambda_\hbar>0$ and $A_2={\tau_\Def}$ if $f\in\mathcal{H}_\infty$ with $\mathrm{supp}(f_\Def)\subset[0,1]$ respectively corresponds to
\begin{equation}\label{SpecialHeisenberg1}
\sum_{\Def\in(-\infty,\infty)}\langle\widehat{\mathcal{Z}}^{\lambda_\Def}_\pm{f_\Def},\widehat{\mathcal{Z}}^{\lambda_\Def}_\pm{f_\Def}\rangle_\Def\geq\Biggr\vert\sum_{\Def\in(-\infty,\infty)}\langle{f_\Def},\mathcal{Z}^{\lambda_\Def}_+f_\Def\rangle_\Def\Biggr\vert^2+\Biggr\vert\sum_{\Def\in(-\infty,\infty)}\langle{f_\Def},\mathcal{Z}^{\lambda_\Def}_-f_\Def\rangle_\Def\Biggr\vert^2
\end{equation}
\end{proposition}
\begin{proof}[Proof] We only use properties of ${\tau_\Def}$ and that $\Ps^\lambda$ is a sum of dilations $\Di_\beta$ with $\beta\geq1$: The vanishing $\langle{f},{\tau_\Def}{f}\rangle_\Def=0$ is a straight forward calculation. It is immediate that $\langle{f},{\tau_\Def}{f}\rangle_\Def=0$ implies for the variance
$\sigma^2_f({\tau_\Def})=\langle{f},\tau_\Def\circ\tau_\Def{f}\rangle_\Def-\langle{f},{\tau_\Def}{f}\rangle_\Def=\langle{f},{f}\rangle_\Def=1$
and \ref{SchroeUncer} reduces to \ref{SchroeUncer2}.

In particular $\mathrm{supp}(f_\Def)\subset[1,\infty]$ or $\subset[0,1]$ implies $\langle{f},{\tau_\Def}{f}\rangle_\Def=0$ and let us now justify the positivity statement \ref{SpecialHeisenberg1}: It is not difficult to verify $\frac{\I}{2}[{\tau_\Def},\widehat{\mathcal{Z}}^{\lambda}_\pm]_{-}=\mp\mathcal{Z}_-^\lambda$ and $\frac{1}{2}[{\tau_\Def},\widehat{\mathcal{Z}}^{\lambda}_\pm]_{+}=\mathcal{Z}_+^\lambda$.
Because of the implication \ref{Support} the inequality \ref{SchroeUncer2} reduces to \ref{SpecialHeisenberg1} respectively and this equations are equivalent because of ${\tau_\Def}\circ\mathcal{Z}^\lambda_\pm\circ{\tau_\Def}=\pm\mathcal{Z}^\lambda_\pm$ and $\mathrm{supp}(f)\subset(1,\infty]\Leftrightarrow\mathrm{supp}({\tau_\Def}{f})\subset[0,1)$.  Analogous if $\langle{f},{\tau_\Def}{f}\rangle_\Def=0$ we find $\sigma^2_f(\mathcal{Z}^\lambda_\pm)\geq\bigr\vert\bigr\langle{f},(\widehat{\mathcal{Z}}^\lambda_+\pm\widehat{\mathcal{Z}}^\lambda_-)f\bigr\rangle\bigr\vert^2$.\end{proof}
Symmetrization of convolution operators \ref{SACon} and the four symmetrizations of the operator $H_\lambda\circ\Ps^\lambda$ defined in \ref{Ad-Hoc} seem similar, but the analogy of \ref{SchroeUncer2} for the corresponding convolution operators seems like a non-trivial question because here we don't have \ref{Support}, it is only immediate to verify the standard  inclusion $\mathrm{supp}(f_1\Conv{f_2})\subset\mathrm{supp}(f_1)\cdot\mathrm{supp}(f_2):=\{r_1\cdot{r_2}\vert{r_i}\in\mathrm{supp}(f_i)\}$.

\subsection{\textit{Involutive construction of the $E\rightarrow-E$ eigenvalue symmetry}}\label{InvoCon}
The easiest method to incorporate the eigenvalue symmetry
$E\rightarrow\pm{E}$
for any self-adjoint operator is by commutation or anti-commutation with ${\tau_\Def}$ respectively:
For $A\in\mathfrak{s}$ we have $[A,{\tau_\Def}]_+,\I[A,{\tau_\Def}]_-\in\mathfrak{s}$ and
$[A,{\tau_\Def}]_\pm\psi_E=E\psi_E\Leftrightarrow[A,{\tau_\Def}]_\pm{\tau_\Def}\psi_E=\pm{E}{\tau_\Def}\psi_E$.
\begin{proposition}\label{EigenSym} $$\left(\mathcal{Z}^\lambda_\mp+\Mu_{V^\mp_m+\overline{V}^\mp_m}+\Co_{V^\mp_c-\overline{V}^\mp_c}\right)\psi_E=E\psi_E\Leftrightarrow\left(\mathcal{Z}^\lambda_\mp+\Mu_{V^\mp_m+\overline{V}^\mp_m}+\Co_{V^\mp_c-\overline{V}^\mp_c}\right){\tau_\Def}\psi_{E}=\mp{E}{\tau_\Def}\psi_{E}$$\end{proposition} 
Proposition \ref{EigenSym} shows that the $\mathcal{Z}^\lambda_-+\Mu_{V^-_m+\overline{V}^-_m}+\Co_{V^-_c-\overline{V}^-_c}$ eigenvalues can only be non-zero if ${\tau_\Def}\psi_E\neq\pm\psi_E$, but for the Hamiltonian $\mathcal{Z}^\lambda_++\Mu_{V^+_m+\overline{V}^+_m}+\Co_{V^+_c-\overline{V}^+_c}$ all the speculative eigenstates have to decompose in two eigenstates $\psi^+_E\in\mathcal{H}^+_\Def$ and $\psi^-_E\in\mathcal{H}^-_\Def$ in the two eigenspaces of $\tau_\Def$.

\section{\textit{Proof II, a closer look at the symmetrization}}\label{proof2}
The second promised proof that $\widehat{\mathcal{Z}}^\lambda_\pm\in\mathfrak{s}$ is more conceptional, just relies on standard procedures and we want to emphasise this way because some other interesting maps arise: It is well-known that for $A,B\in\mathfrak{s}$ also real linear combinations and $\I[A,B]_{-}$ as well as $[A,B]_{+}$ are self-adjoint, the two binary operations $\I[\cdot,\cdot]_-$ and $[\cdot,\cdot]_+$ endow $\in\mathfrak{s}$ with the structure of a \textit{\textbf{Lie-algebra}} and \textit{\textbf{Jordan-algebra}} respectively. In words we can deduce by real linear combinations of
\begin{align}\label{commutator1}
&2\I[\mathcal{Z}^\lambda_{\pm},{\tau_\Def}]_{-}=-(1\mp1)[H_\lambda\circ\Ps^\lambda,{\tau_\Def}]_{-}\quad\text{and}\quad2[\mathcal{Z}^\lambda_{\pm},{\tau_\Def}]_{+}=(1\pm1)[H_\lambda\circ\Ps^\lambda,{\tau_\Def}]_{+}
\end{align}
that $\widehat{\mathcal{Z}}^\lambda_\pm\in\mathfrak{s}$. Notice also $H_\lambda\circ\Ps^\lambda=\widehat{\mathcal{Z}}^\lambda_+\circ{\tau_\Def}={\tau_\Def}\circ\widehat{\mathcal{Z}}^\lambda_-$ and we have in some sense dual to \ref{commutator1}
\begin{align}\label{commutator2}
&\I[[{H}_\lambda\circ\Ps^\lambda,{\tau_\Def}]_{-},{\tau_\Def}]_{\pm}=2(1\mp1)\mathcal{Z}^\lambda_-\quad\text{and}\quad[[{H}_\lambda\circ\Ps^\lambda,{\tau_\Def}]_{+},{\tau_\Def}]_{\pm}=2(1\pm1)\mathcal{Z}^\lambda_+
\end{align}
The standard construction scheme of $\widehat{\mathcal{Z}}^\lambda_\pm$ that lurks in the background is now represented by
\begin{equation}\label{Map}\xymatrix{\mathfrak{g}\ar[rr]^{\hspace{-0.45cm}\left(\mathrm{p},\partial\right)}&&\mathfrak{s}^{\times2}\ar[rr]^{\hspace{+0.05cm}\left([{\tau_\Def},\cdot]_{+},{\I[{\tau_\Def},\cdot]_{-}}\right)}&&\mathfrak{s}^{\times2}\ar[rr]^{\hspace{+0.30cm}{\mathrm{pr}_1\pm\mathrm{pr}_2}}&&\mathfrak{s}}
\end{equation}
We investigate some computations on the moral of the hidden procedure \ref{Map}:

Notice, if we set $[A,B]^\ast_\pm:=\left(A\circ{B}\pm{B}^\ast\circ{A}^\ast\right)$ we can obviously define two binary operations
\begin{equation}
[\cdot,\cdot]_+/2:=\mathrm{p}\circ,\I[\cdot,\cdot]_-/2:=\partial\circ:\mathfrak{g}^{\otimes2}\rightarrow\mathfrak{s}(\Def)
\end{equation}
where the tensor product is over $\mathbb{R}$ and $\circ:\mathfrak{g}^{\otimes2}\rightarrow\mathfrak{s}(\Def)$ maps $(A,B)$ to the composition ${A}\circ{B}$.

Explicit the two maps $\mathfrak{g}\rightarrow\mathfrak{s}(\Def)$ sending $H_\lambda\circ\Ps^\lambda$ to $\widehat{\mathcal{Z}}^\lambda_\pm(\Def)$ described in the sequence \ref{Map} act by
$A\rightarrow({\tau_\Def}\circ{A}+A^\ast\circ{\tau_\Def})/2:=\frac{1}{2}[{\tau_\Def},\cdot]_+^\ast(A)$
and the ${\tau_\Def}$ conjugated formula
$A\rightarrow({\tau_\Def}\circ{A}^\ast+A\circ{\tau_\Def})/2$.
Hence restricted to $\mathfrak{s}(\Def)$ this maps are just $[{\tau_\Def},\cdot]_+/2$ and $\mathrm{Con}_{\tau_\Def}\circ[{\tau_\Def},\cdot]_+/2$. Notice that there is also a negative version of this two maps, namely
$A\rightarrow\I({\tau_\Def}\circ{A}-A^\ast\circ{\tau_\Def})/2:=\frac{\I}{2}[{\tau_\Def},\cdot]_-^\ast(A)$
and minus the ${\tau_\Def}$ conjugated formula
$A\rightarrow\I({\tau_\Def}\circ{A}^\ast-A\circ{\tau_\Def})/2$,
restricted to $\mathfrak{s}(\Def)$
this negative maps are $\I[{\tau_\Def},\cdot]_-/2$ and its ${\tau_\Def}$ conjugate. The negative maps can essentially be recovered from the positive maps by substituting $A\in\mathfrak{g}\rightarrow\I{A}\in\mathfrak{g}$ and vice versa by a small computation.
 \begin{lemma}\label{IterPos}
 Let $\mathcal{I}\in\mathfrak{g}$ be an involution. The following iteration formulas hold:
 $$\left(\begin{cases}{1}/{2}\\{\I}/{2}\end{cases}\hspace{-0.3cm}[\mathcal{I},\cdot]_\pm^\ast\right)^{n+1}=\begin{cases}(\pm1)^{n/2}\left({\id\pm\mathrm{Con}^\ast_\mathcal{I}}\right)\bigr/{2}\hspace{0.5cm}\text{if $n$ is even}\\{(\pm1)^{\frac{n+1}{2}}\begin{cases}{1}/{2}\\{\I}/{2}\end{cases}\hspace{-0.3cm}[\mathcal{I},\cdot]^\ast_\pm}\hspace{0.62cm}\text{if $n$ is odd}\end{cases}\circ\begin{cases}{1}/{2}\\{\I}/{2}\end{cases}\hspace{-0.3cm}[\mathcal{I},\cdot]_\pm^\ast$$\end{lemma}
\begin{proof}[Proof]
For an involution $\mathcal{I}\in\mathfrak{g}$ the equations $\frac{1}{2}[\mathcal{I},\cdot]^\ast_{+}\circ\frac{\I}{2}[\mathcal{I},\cdot]^\ast_{-}=\frac{\id+\mathrm{Con}^\ast_\mathcal{I}}{2}\circ\partial$, $\frac{\I}{2}[\mathcal{I},\cdot]^\ast_{-}\circ\frac{1}{2}[\mathcal{I},\cdot]^\ast_{+}=\frac{\id-\mathrm{Con}^\ast_\mathcal{I}}{2}\circ\partial$ and 
 $(\id\pm\mathrm{Con}_\mathcal{I}^\ast)^2=2(\id\pm\mathrm{Con}_\mathcal{I}^\ast)$, $[\mathcal{I},\cdot]^\ast_{\pm}\circ\mathrm{Con}^\ast_\mathcal{I}=\pm[\mathcal{I},\cdot]^\ast_{\pm}\circ^\ast=\mathrm{Con}^\ast_{\mathcal{I}}\circ[\mathcal{I},\cdot]^\ast_{\pm}$ hold.\end{proof}
\begin{proposition}\label{HoPro}
Let $\mathcal{I}\in\mathfrak{g}$ be an involution. If $\mathcal{I}\circ(A-A')=(A-A')^\ast\circ\mathcal{I}^\ast$ then the negative maps $\I[\mathcal{I},\cdot]^\ast_-/2,\mathrm{Con}^\ast_\mathcal{I}\circ\I[\mathcal{I},\cdot]^\ast_-/2$ will send $A$ and $A'$ to the same element of $\mathfrak{s}(\Def)$, {\em i.e.} the space
$$\mathrm{C}^\ast_{\mathfrak{g}}(\mathcal{I}):=\{A\in\mathfrak{g}\vert\mathcal{I}\circ{A}=A^\ast\circ\mathcal{I}^\ast\}=\{\left({A}+\mathrm{Con}^\ast_\mathcal{I}(A^\ast)\right)/2,A\in\mathfrak{g}\}\subset\mathfrak{g}$$
called the adjoint centralizer of $\mathcal{I}$ is the kernel of the negative maps $\frac{\I}{2}[\mathcal{I},\cdot]_-^\ast$ and its $\mathcal{I}$ adjoint conjugate. Analog $\I\mathrm{C}^\ast_{\mathfrak{g}}(\mathcal{I})$
is the kernel of the two positive maps $\frac{1}{2}[\mathcal{I},\cdot]_+^\ast$ and $\mathrm{Con}^\ast_\mathcal{I}\circ\frac{1}{2}[\mathcal{I},\cdot]_+^\ast$.
\end{proposition}
We can restate \ref{HoPro} saying that the two negative maps are well-defined injective operators $\mathfrak{g}/\mathrm{C}^\ast_{\mathfrak{g}}(\mathcal{I})\rightarrow\mathfrak{s}(\Def)$ and analogous the positive maps descend to well-defined maps $\mathfrak{g}/\I\mathrm{C}^\ast_{\mathfrak{g}}(\mathcal{I})\rightarrow\mathfrak{s}(\Def)$.

The projection $A\rightarrow\left(A+\mathrm{Con}_{\tau_\Def}(A^\ast)\right)/2\in\mathrm{C}^\ast_{\mathfrak{g}}({\tau_\Def})$ does not respect $\frac{\I}{2}[\cdot,\cdot]^\ast_-$. $\mathrm{C}^\ast_{\mathfrak{g}}(\mathcal{I})$ is not closed under adjunction $^{\ast(\Def)}$ or composition, but still a sub Lie- or Jordan-algebra with respect to the usual Lie bracket $\frac{\I}{2}[\cdot,\cdot]_-$ or $\frac{1}{2}[\cdot,\cdot]_+$.

The two inclusions $\tau_\Def\circ{A},A\circ\tau_\Def\in\mathrm{C}^\ast_{\mathfrak{g}}(\tau_\Def)\;\forall{A}\in\mathfrak{s}(\Def)$ and the complementary inclusions $\tau_\Def\circ{A},A\circ\tau_\Def\in\mathfrak{s}(\Def)\;\forall{A}\in\mathrm{C}^\ast_{\mathfrak{g}}(\tau_\Def)$ allow to identify the space $\mathfrak{s}(\Def)$ in different ways with $\mathrm{C}^\ast_{\mathfrak{g}}(\tau_\Def)$.

It is well-known that the eigenvalues of operators contained in $\mathfrak{s}(\Def)$ are real but this is not unconditionally true for operators contained in $\mathrm{C}^\ast_{\mathfrak{g}}(\tau_\Def)$ but almost:
\begin{proposition}
For $A\in\mathrm{C}^\ast_{\mathfrak{g}}(\tau_\Def)$ the eigenvalue equation $A\psi_E=E\psi_E$ implies $E\in\mathbb{R}$ or the vanishing $0=\langle\psi_E,\tau_\Def\psi_E\rangle_\Def\Leftrightarrow0=\langle\psi^+_E,\psi^+_E\rangle_\Def-\langle\psi^-_E,\psi^-_E\rangle_\Def$ with $\psi^\pm_E\in\mathcal{H}_\Def^\pm$ according to \ref{IotaProjection}.
\end{proposition}\begin{proof}[Proof]
In fact the proof only assumes $\langle\alpha{f},\beta{g}\rangle_\Def=\overline{\alpha}\cdot{\beta}\langle{f},{g}\rangle_\Def$ and we do not need that $\langle\cdot,\cdot\rangle_\Def$ is positive: Assume $A\psi_E=E\psi_E$, $A\psi_{E'}=E'\psi_{E'}$ for $A\in\mathrm{C}^\ast_{\mathfrak{g}}(\tau_\Def)$ and we have
$$\overline{E}\langle\psi_E,\tau_\Def\psi_{E'}\rangle_\Def\hspace{-0.04cm}=\hspace{-0.04cm}\langle{A}\psi_E,\tau_\Def\psi_{E'}\rangle_\Def\hspace{-0.04cm}=\hspace{-0.04cm}\langle\tau_\Def\hspace{-0.027cm}\circ\hspace{-0.027cm}{A}^{\ast(\Def)}\hspace{-0.027cm}\circ\hspace{-0.027cm}\tau_\Def\psi_E,\tau_\Def\psi_{E'}\rangle_\Def\hspace{-0.04cm}=\hspace{-0.04cm}\langle\tau_\Def\psi_E,A\psi_{E'}\rangle_\Def\hspace{-0.04cm}=\hspace{-0.04cm}E'\langle\psi_E,\tau_\Def\psi_{E'}\rangle_\Def$$\end{proof}
\begin{proposition}
For an involution $\mathcal{I}$ the maps $\frac{1}{2}[\mathcal{I},\cdot]_+^\ast,\mathrm{Con}^\ast_\mathcal{I}\circ\frac{1}{2}[\mathcal{I},\cdot]_+^\ast=\frac{1}{2}[\mathcal{I},\cdot]_+^\ast\circ^\ast$ and $\frac{\I}{2}[\mathcal{I},\cdot]^\ast_-,\mathrm{Con}^\ast_\mathcal{I}\hspace{-0.02cm}\circ\hspace{-0.03cm}\frac{\I}{2}[\mathcal{I},\cdot]^\ast_-\hspace{-0.1cm}=\hspace{-0.1cm}\frac{\I}{2}[\mathcal{I},\cdot]^\ast_-\hspace{-0.04cm}\circ\hspace{-0.03cm}^\ast$ commute with adjoint conjugation by elements of the centralizer
$$\mathrm{C}_{\mathfrak{g}}(\mathcal{I}):=\{A\in\mathfrak{g}\vert\mathcal{I}\circ{A}=A\circ\mathcal{I}\}=\{\left({A}+\mathrm{Con}_\mathcal{I}(A)\right)/2,{A}\in\mathfrak{g}\}\subset\mathfrak{g}$$
For example
$\mathrm{Con}^\ast_{A}\circ[\mathcal{I},\cdot]_+^\ast=[\mathcal{I},\cdot]_+^\ast\circ\mathrm{Con}^\ast_{A}$ and also $[\mathrm{Con}_{\mathcal{I}}^\ast\circ^\ast,\mathrm{Con}^\ast_A]_-=0$ hold $\forall{A}\in\mathrm{C}_{\mathfrak{g}}(\mathcal{I})$.
\end{proposition}
The vector space $\mathrm{C}_{\mathfrak{g}}(\mathcal{I})$ is closed under adjunction $^{\ast(\Def)}$ if $\mathcal{I}^\ast=\mathcal{I}$ and also under composition and the natural Lie-bracket, we have for $A,B\in\mathrm{C}_{\mathfrak{g}}(\mathcal{I})$ clearly $A\circ{B}\circ\mathcal{I}=\mathcal{I}\circ{A}\circ{B}$.

We have a self-adjoint intersection
$\mathrm{C}^\ast_{\mathfrak{g}}({\tau_\Def})\cap\mathrm{C}_{\mathfrak{g}}({\tau_\Def})\subset\mathfrak{s}(\Def)$. The spaces $\mathrm{C}^\ast_{\mathfrak{g}}({\tau_\Def})=\left(\mathrm{C}^\ast_{\mathfrak{g}}({\tau_\Def})\right)^\ast$ and $\mathrm{C}_{\mathfrak{g}}({\tau_\Def})=\left(\mathrm{C}_{\mathfrak{g}}({\tau_\Def})\right)^\ast$ both contain for example ${\mathcal{Z}}^\lambda_+$ and $\Mu_{V(t)}$ with real potential conjugate symmetric under $t\rightarrow1/t$, hence the friendly \textit{\textbf{Brujin-Newman-P\'olya}} operators $\Mu_{e^{-\rho\ln^2(t)}}$ with $\rho\geq0$ earlier mentioned  while the operators $H_\lambda\circ\Ps^\lambda$ and $\I{\mathcal{Z}}^\lambda_-$ are only contained in $\mathrm{C}^\ast_{\mathfrak{g}}({\tau_\Def})$.
\begin{proposition} For every $V\in\mathcal{H}$ the maps $\Co_{(V+\overline{V})/2}$ are elements of the adjoint centralizer $\mathrm{C}^\ast_{\mathfrak{g}}({\tau_\Def})$ of the hermitian involution ${\tau_\Def}$ and $\Co_{(V+\tau_\Def{V})/2}$ are elements of the centralizer $\mathrm{C}_{\mathfrak{g}}({\tau_\Def})$.\end{proposition}

If $\mathcal{I}\in\mathfrak{g}$ is an involution we have $\left(\mathrm{Con}^\ast_{\mathcal{I}}\right)^2=\id$ and we can analogous to \ref{IotaProjection} also split $A=A^++A^-\in\mathfrak{g}=\mathfrak{g}_\mathcal{I}^+\oplus\mathfrak{g}^-_\mathcal{I}$ where $A^\pm$ are in the two eigenspaces $\mathfrak{g}_\mathcal{I}^\pm$ of $\mathrm{Con}^\ast_\mathcal{I}$, for instance we set
$A^\pm:=\left(A\pm\mathcal{I}\circ{A}\circ\mathcal{I}^\ast\right)/2\in\mathfrak{g}_\mathcal{I}^\pm$.
We have $\mathrm{Con}^\ast_{\tau_\Def}\widehat{\mathcal{Z}}^\lambda_\pm=\widehat{\mathcal{Z}}^\lambda_\mp$ {\em i.e.} $\widehat{\mathcal{Z}}^\lambda_\pm$ are in the same conjugacy class while
${\tau_\Def}\circ\mathcal{Z}^\lambda_\pm\circ{\tau_\Def}=\pm\mathcal{Z}^\lambda_\pm$
{\em i.e.} $\mathcal{Z}^\lambda_\pm$ are in the respective eigenspaces $\mathfrak{g}_{\tau_\Def}^\pm$ of $\mathrm{Con}^\ast_{\tau_\Def}$.

It is also immediate to distinguish two quite familiar operators: We showed the adjunction
\begin{equation}\label{CohTheo}
\frac{\;^\ast\pm\mathrm{Con}_{\tau_\Def}}{2}\begin{cases}\I{H}_\lambda\circ\Ps^\lambda\\H_\lambda\circ\Ps^\lambda\end{cases}=\frac{\mathrm{p}+\I\partial\pm\mathrm{Con}_{\tau_\Def}}{2}\begin{cases}\I{H}_\lambda\circ\Ps^\lambda\\H_\lambda\circ\Ps^\lambda\end{cases}=0
\end{equation}
and rewriting $\mathrm{Con}^\ast_{\tau_\Def}=\mathrm{Con}_{\tau_\Def}$ in \ref{CohTheo} is more practical: Because$\;^\ast$ is an involution and commutes with $\mathrm{Con}_B^\ast$ we have for all $B\in\mathfrak{g}$ the formula $\left(\frac{\;^\ast\pm\mathrm{Con}^\ast_B}{2}\right)^2=\frac{\id+\mathrm{Con}^\ast_{B\circ{B}}\pm2\mathrm{Con}^\ast_{B}\circ^\ast}{4}$.
\begin{lemma}
Let the operator $\mathcal{I}\in\mathfrak{g}$ be an involution {\em i.e.} $\mathcal{I}\circ{\mathcal{I}}=\id$. We have
$$\left(\frac{\;^\ast\pm\mathrm{Con}^\ast_\mathcal{I}}{2}\right)^{n}=\begin{cases}\frac{\;^\ast\pm\mathrm{Con}_\mathcal{I}^\ast}{2}\;\quad\quad\text{if $n$ is odd}\\\frac{\id\pm\mathrm{Con}_\mathcal{I}^\ast\circ^\ast}{2}\quad\text{if $n$ is even}\end{cases}$$
and the operators $\left({\id\pm\mathrm{Con}_\mathcal{I}^\ast\circ^\ast}\right)/{2}$ are projectors. By flipping signs we have in fact two projectors
$$\left(\id-\frac{\id\pm\mathrm{Con}_\mathcal{I}^\ast\circ^\ast}{2}\right)=\frac{\id\mp\mathrm{Con}_\mathcal{I}^\ast\circ^\ast}{2}:\mathfrak{g}\rightarrow\mathrm{ker}\frac{\;^\ast\pm\mathrm{Con}_\mathcal{I}^\ast}{2}=\frac{\I(1\pm1)+1\mp1}{2}\mathrm{C}^\ast_{\mathfrak{g}}(\mathcal{I})$$
onto the solutions $S$ of the respective equations $\frac{\;^\ast\pm\mathrm{Con}^\ast_\mathcal{I}}{2}S=0$.

If $\frac{\;^\ast-\mathrm{Con}^\ast_\mathcal{I}}{2}S=0$ then $\mathcal{I}\circ{S},S\circ{\mathcal{I}}^\ast\in\mathfrak{s}$ and if $\frac{\;^\ast+\mathrm{Con}^\ast_\mathcal{I}}{2}S=0$ then we have $\I{\mathcal{I}}\circ{S},\I{S}\circ{\mathcal{I}}^\ast\in\mathfrak{s}$.
\end{lemma}
\subsection{\textit{The artificial $\Def$-dependence of the symmetrizations}}\label{Hdependence}

Consider the category with one object $\mathcal{H}$ and morphisms the set of adjoinable operators $\mathfrak{g}$.
The map $\Mu_{t^{\Def-\Def'}}:\mathcal{H}\rightarrow\mathcal{H}$ defines an isomorphism $^{\ast(\Def)}\rightarrow^{\ast(\Def')}$ between the two fully faithful contravariant functors $^{\ast(\Def)}$ and $^{\ast(\Def)'}$ because 
$\;^{\ast(\Def')}={\mathrm{Con}_{\Mu_{t^{\Def-\Def'}}}}\circ^{\ast(\Def)}=^{\ast(\Def)}\circ\hspace{0.08cm}{\mathrm{Con}_{\Mu_{t^{\Def'-\Def}}}}$.

The symmetrization procedures $\mathrm{p}(\Def)$ and $\partial(\Def)$ are not compatible with composition, we have 
\begin{align}
&\begin{cases}\mathrm{p}\\\partial\end{cases}\hspace{-0.4cm}(\Def')=\begin{cases}\frac{1}{2}({1-\mathrm{Con}_{\Mu_{t^{\Def-\Def'}}}})\\\frac{\I}{2}({1-\mathrm{Con}_{\Mu_{t^{\Def-\Def'}}}})\end{cases}+\mathrm{Con}_{\Mu_{t^{\Def-\Def'}}}\circ\begin{cases}\mathrm{p}\\\partial\end{cases}\hspace{-0.4cm}(\Def)\nonumber
\end{align}
hence restrictions $(\mathrm{p},\partial)(\Def')=\mathrm{Con}^{\times2}_{\Mu_{t^{\Def-\Def'}}}\circ(\mathrm{p},\partial)(\Def):\mathrm{C}^{\times2}_{\mathfrak{g}}\left({\Mu_{t^{\Def-\Def'}}}\right)\rightarrow\mathfrak{s}^{\times2}(\Def')$.

The space of multiplication operators $\Mu_{V(t)}$ as defined in \ref{MuLti} is a subspace of $\mathrm{C}_{\mathfrak{g}}\left({\Mu_{t^{\Def-\Def'}}}\right)$.
\begin{proposition}\label{Hindependence}
The involutive contravariant functors
$\mathrm{Con}_{\tau_{\Def}}\circ^{\ast(\Def)}=\hspace{0.1cm}^{\ast(\Def)}\circ\hspace{0.1cm}\mathrm{Con}_{\tau_{\Def}}$
do not depend on the deformation parameter $\Def$.
\end{proposition}
\begin{proof}[Proof] First
$\tau_{\Def'}={\mathrm{Con}_{\Mu_{t^{\frac{\Def-\Def'}{2}}}}}(\tau_{\Def})$ is immediate.
Again the two covariant functors $\mathrm{Con}_{\tau_{\Def}}$ and $\mathrm{Con}_{\tau_{\Def'}}$ are connected by the natural transformation $\Mu_{t^{\Def-\Def'}}:\mathcal{H}\rightarrow\mathcal{H}$, we have $\mathrm{Con}_{\tau_{\Def'}}={\mathrm{Con}_{\Mu_{t^{\Def-\Def'}}}}\circ\mathrm{Con}_{\tau_{\Def}}\nonumber=\mathrm{Con}_{\tau_{\Def}}\circ{\mathrm{Con}_{\Mu_{t^{\Def'-\Def}}}}$.
The substitution
$\tau_\Def\circ\Mu_{t^\lambda}=\Mu_{t^{-\lambda}}\circ\tau_\Def$
implies for the composed contravariant functors $\mathrm{Con}_{\tau_{\Def'}}\circ^{\ast(\Def')}=\mathrm{Con}_{\tau_{\Def}}\circ^{\ast(\Def)}$ and the commutative diagram
\vspace{-0.2cm}\begin{equation}\label{TC}
\hspace{-0.cm}\xymatrix{\mathcal{H}\ar@<1ex>[d]^{{\Mu_{t^{\Def-\Def'}}}}&&&\mathcal{H}\ar[rrr]^{\mathrm{Con}_{\tau_{\Def}}A}\ar[lll]_{\quad{A}^{\ast(\Def)}}\ar@<1ex>[d]^{{\Mu_{t^{\Def-\Def'}}}}&&& \mathcal{H}\ar@<1ex>[d]^{{\Mu_{t^{\Def-\Def'}}}}\\ \mathcal{H}\ar@<1ex>[u]^{{\Mu_{t^{\Def'-\Def}}}}&&&
\mathcal{H}\ar@<1ex>[u]^{{\Mu_{t^{\Def'-\Def}}}}\ar[rrr]_{\mathrm{Con}_{\tau_{\Def'}}A}\ar[lll]^{\quad{A}^{\ast(\Def')}}&&& \mathcal{H}\ar@<1ex>[u]^{{\Mu_{t^{\Def'-\Def}}}}}
\vspace{-0.2cm}\end{equation}

The following commutative diagram \ref{TC2} summarizes some considerations, the dashed line represents the contravariant $\Def$ independent functor of \ref{Hindependence}:
\vspace{-0.2cm}\begin{equation}\label{TC2}\xymatrix{&&&&\mathfrak{g}\ar[rrrd]^{\mathrm{Con}_{\tau_{\Def}}}\ar[dllll]_{{\ast(\Def)}}\ar@<1ex>[ddl]^{\mathrm{Con}_{\Mu_{t^{\Def-\Def'}}}}&&&\\ \mathfrak{g}\ar[rrrru]\ar[rrrd]\ar@{.>}@/^5.2pc/[rrrrrrr]\ar@{.>}@/^5.2pc/[rrrrrrr]^{\mathrm{Con}_{\tau_{\Def}}\circ^{\ast(\Def)}=\hspace{0.1cm}\mathrm{Con}_{\tau_{\Def'}}\circ^{\ast(\Def')}}&&&&&&&\mathfrak{g}\ar[lllu]\ar[lllld]\ar@{.>}@/_5.2pc/[lllllll]\ar@{.>}@/_5.2pc/[lllllll]\\&&&
 \mathfrak{g}\ar@<1ex>[uur]^{\mathrm{Con}_{\Mu_{t^{\Def'-\Def}}}}\ar[rrrru]_{\mathrm{Con}_{\tau_{\Def'}}}\ar[lllu]^{^{\ast(\Def')}}&&&& }\vspace{-0.8cm}\end{equation}\end{proof}

We can rewrite the $\Def$-dependence of the respective involutions by
$$\tau_{\Def'}=\Mu_{t^{\Def-\Def'}}\circ\tau_{\Def}=\tau_{\Def}\circ\Mu_{t^{\Def'-\Def}}$$
and consider the $\Def$ dependence of the other two symmetrization procedures described in \ref{Ad-Hoc}, again the natural transformation $^{\ast(\Def)}\rightarrow^{\ast(\Def')}$ given explicit by $\Mu_{t^{\Def-\Def'}}:\mathcal{H}\rightarrow\mathcal{H}$ appears, but as a composition acting on the space of adjoinable operators $\mathfrak{g}$:
\begin{lemma} We have the symmetrization identities 
\begin{align*}
&\begin{cases}\frac{1}{2}[\tau_{\Def'},\cdot]_+^{\ast(\Def')}\\\frac{\I}{2}[\tau_{\Def'},\cdot]^{\ast(\Def')}_-\end{cases}\hspace{-0.4cm}=\Mu_{t^{\Def-\Def'}}\circ\begin{cases}\frac{1}{2}[\tau_{\Def},\cdot]_+^{\ast(\Def)}\\\frac{\I}{2}[\tau_{\Def},\cdot]^{\ast(\Def)}_-\end{cases}\;\text{and}\quad\begin{cases}\frac{1}{2}[\tau_{\Def'},\cdot]_+^{\ast(\Def')}\\\frac{\I}{2}[\tau_{\Def'},\cdot]^{\ast(\Def')}_-\end{cases}\hspace{-0.4cm}\circ^{\ast(\Def')}=\begin{cases}\frac{1}{2}[\tau_{\Def},\cdot]_+^{\ast(\Def)}\\\frac{\I}{2}[\tau_{\Def},\cdot]^{\ast(\Def)}_-\end{cases}\hspace{-0.4cm}\circ^{\ast(\Def)}\circ\Mu_{t^{\Def'-\Def}}
\end{align*}
and can express this for example in the commutative diagram
\begin{equation}\label{TC}\hspace{0.5cm}\xymatrix{\mathfrak{s}(\Def)\ar@/^1.0pc/[rrrrrr]\ar@<1ex>[dd]^{\circ{\Mu_{t^{\Def'-\Def}}}}&&&&&& \mathfrak{s}(\Def)\ar@<1ex>[dd]^{{\Mu_{t^{\Def-\Def'}}}\circ}\ar@/_1.0pc/[llllll]^{\mathrm{Con}_{\tau_{\Def}}}\\&&&\mathfrak{g}/\left(\I\mathrm{C}^\ast_{\mathfrak{g}}(\tau_{\Def})\cap\I\mathrm{C}^\ast_{\mathfrak{g}}(\tau_{\Def'})\right)\ar[urrr]^{\frac{1}{2}[\tau_{\Def},\cdot]_+^{\ast(\Def)}}\ar[ulll]_{\quad\frac{1}{2}[\tau_{\Def},\cdot]_+^{\ast(\Def)}\circ^{\ast(\Def)}}\ar[drrr]_{\frac{1}{2}[\tau_{\Def'},\cdot]_+^{\ast(\Def')}}\ar[dlll]^{\quad\frac{1}{2}[\tau_{\Def'},\cdot]_+^{\ast(\Def')}\circ^{\ast(\Def')}}&&&&&\\\mathfrak{s}(\Def')\ar@/_1.0pc/[rrrrrr]\ar@<1ex>[uu]^{\circ{\Mu_{t^{\Def-\Def'}}}}&&&&&& \mathfrak{s}(\Def')\ar@/^1.0pc/[llllll]_{\mathrm{Con}_{\tau_{\Def'}}}\ar@<1ex>[uu]^{{\Mu_{t^{\Def'-\Def}}}\circ}}
\end{equation}
By the propositions \ref{HoPro} and \ref{Hindependence} we have in fact $\I\mathrm{C}^\ast_{\mathfrak{g}}(\tau_{\Def})=\I\mathrm{C}^\ast_{\mathfrak{g}}(\tau_{\Def'})=\I\mathrm{C}^\ast_{\mathfrak{g}}(\tau_{\Def})\cap\I\mathrm{C}^\ast_{\mathfrak{g}}(\tau_{\Def'})$.\end{lemma}
\section{Certain quotient spaces and cohomologies}\label{QuoEigen}
We proceed quite in analogy with the absorption spectrum of Connes but instead of quotient out the image of the subspace defined by the condition \ref{AltCon} under the operator $\Ps^\lambda$ we quotient out the images of \ref{DH}: Let $\mathcal{H}^\lambda_\pm(\Def)$, $\widehat{\mathcal{H}}^\lambda_\pm(\Def)$ be the quotients of $\mathcal{H}$ by the respective image of the four operators $\mathcal{Z}^\lambda_\pm(\Def)$, $\widehat{\mathcal{Z}}^\lambda_\pm(\Def)$, in formulas $\mathcal{H}^\lambda_\pm:=\mathcal{H}/{\mathcal{Z}}^\lambda_\pm$, $\widehat{\mathcal{H}}^\lambda_\pm:=\mathcal{H}/\widehat{\mathcal{Z}}^\lambda_\pm$, {\em i.e.} elements of this space are equivalence classes $\psi^\lambda_f(\Def)$ with $f\in\mathcal{H}$, we have by definition for two representatives 
\begin{equation}\label{Quotient}
\psi^\lambda_f(\Def)=\psi^\lambda_{f'}(\Def)\Leftrightarrow{f}=f'+\mathcal{Z}^\lambda_\pm(\Def){g}\quad\text{and}\quad\psi^\lambda_f(\Def)=\psi^\lambda_{f'}(\Def)\Leftrightarrow{f}=f'+\widehat{\mathcal{Z}}^\lambda_\pm(\Def){g}
\end{equation}
with $g\in\mathcal{H}$. There are technical issues but in consideration of the following discussion \ref{QuoEigen} it seems desirable to further restrict $g$, maybe Hurwitz stability could play a role.

Because we suppose $\hbar\in\mathbb{R}$ and $\lambda\in\mathbb{R}^+$ complex conjugation $\overline{\psi^\lambda_f(\Def)}=\psi^\lambda_{\overline{f}}(\Def)$ is a involution. Also $\tau_\Def$ induces a well-defined involution $\tau_\Def$ on $\mathcal{H}^\lambda_\pm(\Def)$ by the formula
$\tau_\Def\psi^\lambda_f:=\psi_{\tau_\Def{f}}^\lambda$.
 If $\hbar=\frac{1-\lambda}{\lambda}$ the issue seems more delicate in consideration of $\zeta(1-s)=2(2\pi)^{-s}\cos\left(\frac{\pi{s}}{2}\right)\Gamma(s)\zeta(s)$, however $\tau_\Def:\psi^\lambda_f\rightarrow\psi_{\tau_\Def{f}}^\lambda$ is well-defined  for the quotient $\mathcal{H}/\langle\widehat{\mathcal{Z}}^\lambda_\pm(\hbar)\rangle$ by the linear span of $\widehat{\mathcal{Z}}^\lambda_\pm(\hbar)$. 


In fact the previous defined quotients are not only vector spaces, they are commutative rings: Because of \ref{CoAd-Hoc} and \ref{CoCo1} we have induced by the usual convolution \ref{MuCo} a well-defined associative product
$\Conv:\widehat{\mathcal{H}}^\lambda_\pm\left(\Def\right)\times\widehat{\mathcal{H}}^\lambda_\pm\left(\Def\right)$ and $\Conv:\mathcal{H}^\lambda_\pm\left(\Def\right)\times\mathcal{H}^\lambda_\pm\left(\Def\right)\rightarrow\mathcal{H}^\lambda_\pm\left(\Def\right)$
naturally by the same formula
\begin{equation}
\psi^\lambda_{f_1}\Conv\psi^\lambda_{f_2}:=\psi^\lambda_{f_1\Conv{f}_2}
\end{equation}

 The computations in \ref{Std1} show that $H_\alpha$ descend by $H_\alpha\psi^\lambda_f(\Def):=\psi^\lambda_{{H_\alpha}{f}}(\Def)$ to well-defined operators on the four quotient spaces because $H_\alpha$ commutes with $\Ps^\lambda$ and we have \ref{ComRel1}.
 
 \subsubsection{Discussion of eigenvalue equations for the quotients}\label{QuoEigen}
 Let $\mu(n)$ be the M\"obius function defined by $1$ for $n=1$, $0$ if $n$ is not square free and $(-1)^m$ if $n=p_1\cdots{p}_m$ with $p_i\neq{p}_j$. The map $\Ps^\lambda$ admits in some sense an inverse $(\Ps^\lambda)^{-1}:=\sum_{n=1}^\infty\mu(n)\Di_{n^\lambda}$ we refer the reader to \cite{RM} for a precise definition of the domain of this operators, but this operators are not maps $\mathcal{H}\rightarrow\mathcal{H}$, they could produce a singularity at zero. Although $H_\alpha:\mathcal{H}\rightarrow\mathcal{H}$ is injective ($t^{-1/\alpha}$ is not an element of $\mathcal{H}$ hence the kernel of $H_\alpha:\mathcal{H}\rightarrow\mathcal{H}$ trivial) $H_\alpha$ is also not invertible by the natural candidate geometric series, this operation is highly singular.
Hence the existence of non-trivial solutions  of the following eigenvalue equations
\ref{QuotientFunc1} is {\em a priori} not justified: 
\begin{equation}\label{QuotientFunc1}
H_{\alpha}\psi^\lambda_f\left(\hbar\right)=E\psi^\lambda_f\left(\hbar\right)\Leftrightarrow\exists{g}\in\mathcal{H}:H_{\alpha}f-Ef=\begin{cases}\mathcal{Z}^\lambda_\pm\left(\hbar\right)\\\widehat{\mathcal{Z}}^\lambda_\pm\left(\hbar\right)\end{cases}\hspace{-0.3cm}g
\end{equation}
We set $\Def=(1-\lambda)/\lambda$ and $\alpha=2\lambda$ and Mellin transform  of \ref{QuotientFunc1} implies
\begin{align}\label{MeQuo1}&(1-E-2{s}){\mathcal{M}[f]}\left(\frac{s}{\lambda}\right)=\begin{cases}1\\\I\end{cases}\hspace{-0.3cm}\bigr[(1-\id)\cdot\zeta(s)\pm(1-\id)\cdot\zeta(1-s)\bigr]{\mathcal{M}[g]}\left(\frac{s}{\lambda}\right)\\&(1-E-2{s}){\mathcal{M}[f]}\left(\frac{(1\mp1)/2\pm{s}}{\lambda}\right)=(1-\id)\cdot\zeta\bigr((1\mp1)/2\pm{s}\bigr){\mathcal{M}[g]}\left(\frac{(1\pm1)/2\mp{s}}{\lambda}\right)\nonumber\end{align}The l.h.s of the equations \ref{MeQuo1} clearly vanishes at $z=(1-E)/2$ and we yield
\begin{align}\label{VaniMe}&0=\Biggr[(1-\id)\cdot\zeta\left(\frac{1-E}{2}\right)\pm(1-\id)\cdot\zeta\left(\frac{1+E}{2}\right)\Biggr]{\mathcal{M}[g]}\left(\frac{1-{E}}{2\lambda}\right)\\&0=(1-\id)\cdot\zeta\left(\frac{1-E}{2}\right){\mathcal{M}[g]}\left(\frac{1\pm{E}}{2\lambda}\right)\nonumber\end{align}
Suppose $\zeta(z)=0$: By the convolution theorem \ref{ConvolutionTh} and the well-known functional equation $\Xi(s)=\Xi(1-s)$ there exist well-defined sesquilinear maps
$\langle\cdot,\cdot\rangle^\lambda_z:\mathcal{H}^\lambda_\pm\left(\frac{1-\lambda}{\lambda}\right)\times\mathcal{H}^\lambda_\pm\left(\frac{1-\lambda}{\lambda}\right)\rightarrow\mathbb{C}$ and $\langle\cdot,\cdot\rangle^\lambda_z:\mathcal{H}/\langle\widehat{\mathcal{Z}}^\lambda_\pm\rangle\left(\frac{1-\lambda}{\lambda}\right)\times\mathcal{H}/\langle\widehat{\mathcal{Z}}^\lambda_\pm\rangle\left(\frac{1-\lambda}{\lambda}\right)\rightarrow\mathbb{C}$
respectively defined by the same formula
\begin{equation}\Bigr\langle\psi^\lambda_{f_1},\psi^\lambda_{f_2}\Bigr\rangle^\lambda_z:=\mathcal{M}\left[\overline{\tau_{\frac{1-\lambda}{\lambda}}{f_1}}\Conv{f_2}\right](z/\lambda)\end{equation}
We have 
$\langle\psi^\lambda_{\overline{f}_2},\psi^\lambda_{\overline{f}_1}\rangle^\lambda_{1-z}=\Bigr\langle\psi^\lambda_{\tau_\frac{1-\lambda}{\lambda}\overline{f_2}},\psi^\lambda_{\tau_\frac{1-\lambda}{\lambda}\overline{f_1}}\Bigr\rangle^\lambda_z=\langle\psi^\lambda_{f_1},\psi^\lambda_{f_2}\rangle^\lambda_{z}$. The calculations\ref{ComRel1} and \ref{HCONV} imply
\begin{equation}
\Bigr\langle{H_{2\lambda}\psi^\lambda_{f_1}\left(\frac{1-\lambda}{\lambda}\right),\psi^\lambda_{f_2}\left(\frac{1-\lambda}{\lambda}\right)\Bigr\rangle^\lambda_z=-\Bigr\langle\psi^\lambda_{f_1}\left(\frac{1-\lambda}{\lambda}\right),H_{2\lambda}}\psi^\lambda_{f_2}\left(\frac{1-\lambda}{\lambda}\right)\Bigr\rangle^\lambda_{z}
\end{equation}
Hence if $\psi^\lambda_{f_1}$ and $\psi^\lambda_{f_2}$ are eigenstates of $H_{2\lambda}$ with eigenvalues $E_1$ and $E_2$ we have \textit{``orthogonality"}
\begin{equation}\label{NewEigen}
(\overline{E}_1+E_2)\Bigr\langle\psi^\lambda_{f_1}\left(\frac{1-\lambda}{\lambda}\right),\psi^\lambda_{f_2}\left(\frac{1-\lambda}{\lambda}\right)\Bigr\rangle^\lambda_{z}=0
\end{equation}
If $H_{2\lambda}\psi^\lambda_{f}=E\psi^\lambda_{f}$ with $z=\frac{1-E}{2}$ vanishing of \ref{NewEigen} implies $\Re(z)=\frac{1}{2}$ or $0=\langle\psi^\lambda_{f},\psi^\lambda_{f}\rangle^\lambda_z$, but we can not justify the existence of solutions of \ref{QuotientFunc1} with $g$ restricted to have only zeros on the critical line and also $\Conv$ is {\em a priori} not well-defined on such a restricted quotient space. 
 
\subsection{The analogy with Meyer's spectral realizations}\label{QuoEigen2}
We adapt a idea developed in \cite{Co2},\cite{RM}:  The dual spaces  $\left({\mathcal{H}}^\lambda_\pm(\Def)\right)^\ast$, $\bigr(\widehat{\mathcal{H}}^\lambda_\pm(\Def)\bigr)^\ast$ are the spaces of linear functionals $\varphi:\mathcal{H}\rightarrow\mathbb{C}$ that annihilate the image of $\mathcal{Z}^\lambda_\pm(\Def)$, $\widehat{\mathcal{Z}}^\lambda_\pm(\Def)$ respectively, hence spanned by $\varphi^l_{z}$ given by $\psi_f\rightarrow\mathcal{M}\big[\Mu_{\ln^l(t)}f\big](\frac{1+\Def}{2}+\frac{z}{\lambda})$ with $0\leq{l}<k_z\in\mathbb{N}^+$ and where $z$ is a zero of multiplicity $k_z$ of the respective functions $(1-\id)\zeta\left(\frac{(1+\Def)}{2/\lambda}+z\right)\pm(1-\id)\zeta\left(\frac{(1+\Def)}{2/\lambda}-z\right)$, $(1-\id)\zeta\left(\frac{(1+\Def)}{2/\lambda}\hspace{-0.05cm}\pm\hspace{-0.05cm}z\hspace{-0.05cm}\right)$. The cup product  defined by $\big(\varphi_1\cup\varphi_2\big)(\psi^\lambda_{f_1},\psi^\lambda_{f_2}):=\varphi_1(\psi^\lambda_{f_1}(\Def))\varphi_2(\psi^\lambda_{f_2}(\Def))$ satisfies the inequality $\big(\varphi_z\cup\varphi_{\overline{z}}\big)(\psi^\lambda_{f},\psi^\lambda_{\overline{f}})\geq0$. The convolution theorem implies
\begin{align}\label{cup1}\big(\varphi^{l_1}_{z_1}\cup\varphi^{l_2}_{z_2}\big)(\psi^\lambda_{f_2},\psi^\lambda_{f_2})&=\mathcal{M}\Big[\Mu_{t^{\frac{z_1-z_2}{2\lambda}}\ln^{l_1}(t)}f_1\Conv\Mu_{t^{\frac{z_2-z_1}{2\lambda}}\ln^{l_2}(t)}f_2\Big]\left(\frac{1+\Def}{2}+\frac{z_1+z_2}{2\lambda}\right)\end{align}
Hence the transposed operators $^{t}\hspace{-0.05cm}H_{\alpha}$  acting on the respective dual spaces satisfy
\begin{align}\label{cupH}^{t}\hspace{-0.05cm}H_{\alpha}\varphi^{l_1}_{z_1}\cup\varphi^{l_2}_{z_2}-\varphi^{l_1}_{z_1}\cup^{t}\hspace{-0.05cm}H_{\alpha}\varphi^{l_2}_{z_2}=&\frac{\alpha(z_2-z_1)}{\lambda}\varphi^{l_1}_{z_1}\cup\varphi^{l_2}_{z_2}+\alpha\left(l_2\varphi^{l_1}_{z_1}\cup\varphi^{l_2-1}_{z_2}-l_1\varphi^{l_1-1}_{z_1}\cup\varphi^{l_2}_{z_2}\right)\end{align}
\begin{align}\label{cupH}^{t}\hspace{-0.05cm}H_{\alpha}\varphi^{l}_{z}(\psi^\lambda_{f})&=\mathcal{M}\big[\Mu_{\ln^l(t)}H_{\alpha}f\big]\left(\frac{1+\Def}{2}+\frac{z}{\lambda}\right)\nonumber\\&=\mathcal{M}\big[H_{\alpha}\Mu_{\ln^l(t)}f\big]\left(\frac{1+\Def}{2}+\frac{z}{\lambda}\right)-\alpha{l}\mathcal{M}\big[\Mu_{\ln^{l-1}(t)}f\big]\left(\frac{1+\Def}{2}+\frac{z}{\lambda}\right)\end{align}
where we used \ref{cup1} and the Leibniz rule. By integration by parts $\mathcal{M}\big[H_{\alpha}f\big](s)=(1-\alpha{s})\mathcal{M}\big[f\big](s)$ holds for $f\in\mathcal{H}$ and \ref{cupH} implies the formula $\left(^{t}\hspace{-0.05cm}H_{\alpha}-1+\alpha\left(\frac{1+\Def}{2}+\frac{z}{\lambda}\right)\right)\varphi^{l}_{z}=-\alpha{l}\varphi^{l-1}_{z}$.
Hence if we suppose $(1-\id)\zeta\left(\frac{(1+\Def)}{2/\lambda}+z\right)\pm(1-\id)\zeta\left(\frac{(1+\Def)}{2/\lambda}-z\right)=0$ or $0=(1-\id)\zeta\left(\frac{(1+\Def)}{2/\lambda}\hspace{-0.05cm}\pm\hspace{-0.05cm}z\hspace{-0.05cm}\right)$ respectively
the functionals  $\varphi^l_{z}$ with $0\leq{l}<k_z$ are in the kernel of $(^{t}\hspace{-0.05cm}H_{\alpha}-1+\alpha\left(\frac{1+\Def}{2}+\frac{z}{\lambda}\right))^{l+1}$. In particular
$\varphi^{0}_{z}$ is a eigenstate of $H_{\alpha}$ with the eigenvalue $1-\alpha\left(\frac{1+\Def}{2}+\frac{z}{\lambda}\right)$ and the dimension of $\cup_{l\in\mathbb{N}}\ker\left(^{t}\hspace{-0.05cm}H_{\alpha}-1+\alpha\left(\frac{1+\Def}{2}+\frac{z}{\lambda}\right)\right)^{l+1}$ is $k_z$.
\subsection{Some cohomologies}\label{QuoEigen3}
In the following lemmas $(\mathcal{R},\Conv,+)$ denotes a commutative ring, $\tau:\mathcal{R}\rightarrow\mathcal{R}$ an involutive ring homomorphism and $\mathcal{R}^\pm$ denote the eigenspaces of $\tau$ for the eigenvalues $\pm1$.
\begin{lemma}\label{LieBra}
The formulas
$\star^\pm:=\Conv\circ(\id\pm\tau)\otimes(\id\pm\tau)$
define associative products and $f_1\star^\pm{f_2}\in\mathcal{R}^+$. We have a Lie-bracket  
$[\cdot,\cdot]^\ast:=\Conv\circ\bigr(\tau\otimes\id-\id\otimes\tau\bigr)$ and $[f_1,f_2]^\ast\in\mathcal{R}^-$.
$\tau$ is a morphism of $\star^\pm$ and a Lie morphism of $[\cdot,\cdot]^\ast$.
\end{lemma}
\begin{proof}[Proof] $\star^\pm=\Conv\circ(\id\pm\tau)^{\otimes2}$ also is associative if $\tau$ symbolizes a involutive morphism of a non-commutative product $\Conv$.  Observe if $\Conv$ is commutative $f\Conv_{V}g:=f\Conv{g}\Conv{V}$ is also associative. 
\end{proof}
The following lemma provides two linear maps$\mathcal{R}/\mathcal{R}^+\rightarrow\bigr\{\d:\mathcal{R}\hspace{-0.15cm}\xymatrix{\tiny\ar@(dr,ur)}\quad\quad\vert\;\d\;\text{linear},\;\d^2=0\bigr\}$.
\begin{lemma}\label{Differential}
The maps
$\d^+_V:=[(\id-\tau)V,\cdot]^\ast:\mathcal{R}\rightarrow\mathcal{R}$
are inner derivations of $[\cdot,\cdot]^\ast$ and $\d^-_V$ defined by $\d^-_V{f}:=V\star^-{f}$ satisfies $\d^-_V({f\Conv{g}})=\d^-_V{f}\Conv{g}+\tau{f}\Conv\d^-_V{g}$. The equations $\d^\pm_{V_1}\circ\d^\pm_{V_2}=0$ hold $\forall\;{V_1,V_2}\in\mathcal{R}$ and we have the inclusion $\mathcal{R}^\mp\subset\ker(\d_V^\pm)$.
The product rule $\d^\pm_{V_1\Conv{V}_2}f=V_1\Conv\d^\pm_{V_2}{f}+\tau{V_2}\Conv\d^\pm_{V_1}f$ holds, hence $\ker(\d_{V_1}^\pm)\cap\ker(\d_{V_2}^\pm)\subset\ker(\d_{V_1\Conv{V}_2}^\pm)$.
\end{lemma}
\begin{proof}[Proof] More general if $V_1\in\mathcal{R}^-$ and $V_1\Conv{V}_1=-V_2\Conv\tau{V}_2$ we have a square zero map $\d_{V_1,V_2}$ by setting $\d_{V_1,V_2}{f}:=V_1\Conv{f}+V_2\Conv\tau{f}$. In particular we have $\d_{(\id-\tau)V,\pm(\id-\tau)V}=\d^\pm_V$.
\end{proof}
As usual we can induce a Lie-bracket $[\cdot,\cdot]^\ast$ on the cohomology $H(\mathcal{R},d^+_V):=\mathrm{ker}(\d^+_V)/\mathrm{im}(\d^+_V)$ by
$[\lfloor{f_1}\rfloor_V^+,\lfloor{f_2}\rfloor_V^+]^\ast:=\lfloor[f_1,f_2]^\ast\rfloor_V^+$
for representatives $\lfloor{f_i}\rfloor_V^+\in{H}(\mathcal{R},d^+_V)$. Also  $\Conv$ induces a product on $H(\mathcal{R},d^-_V):=\mathrm{ker}(\d^-_V)/\mathrm{im}(\d^-_V)$ by
$\lfloor{f_1}\rfloor_V^-\Conv\lfloor{f_2}\rfloor_V^-:=\lfloor{f}_1\Conv{f}_2\rfloor_V^-$ for $\lfloor{f_i}\rfloor_V^-\in{H}(\mathcal{R},d^-_V)$. We have a involution on the cohomologies by $\tau\lfloor{f}\rfloor_V^\pm:=\lfloor\tau{f}\rfloor_V^\pm$.

\begin{proposition}\label{dH} For $\mathcal{R}=\mathcal{H}^\lambda_\pm(\Def)$ and $\mathcal{R}=\mathcal{H}/\langle\widehat{\mathcal{Z}}^\lambda_\pm\rangle$ respectively the operator  $H_{\frac{2}{1+\Def}}$ defines a map $H(\mathcal{R},\d^\pm_V)\rightarrow{H}(\mathcal{R},\d^\mp_V)$ and $\Delta_{\frac{2}{1+\Def}}\lfloor{\psi_f}\rfloor_V:=\lfloor{\psi_{\Delta_{\frac{2}{1+\Def}}f}}\rfloor_V$ descends to $H(\mathcal{R},\d^\pm_V)$.
\end{proposition}
\begin{proof}[Proof] We have shown \ref{ComRel1} that $H_{{2}/({1+\Def})}$ anti-commutes with $\tau_\Def$, $\widehat{\mathcal{Z}}^\lambda_\pm(\Def)$ and commutes with $\mathcal{Z}^\lambda_\pm(\Def)$ and we have the identity ${\Delta}_\alpha=H_\alpha^2-1=(2\alpha+\alpha^2){t}\partial_{t}+\alpha^2t^2\partial_{t}^2$.
\end{proof}
For the usual convolution \ref{MuCo} we have the adjunction $(\d^\pm_V)^\ast=-\d^\mp_{\overline{V}}$ with respect to the inner product \ref{SP}, because $\d^\pm_V=\Co_{V-\tau{V}}\circ(\id\pm\tau)$ and \ref{CoAd}, \ref{CoCo1} implies $(\d^\pm_V)^\ast=(\id\pm\tau)\circ\Co_{\overline{V}-\tau\overline{V}}=\d^\mp_{\overline{V}}$. 

Observe if $V-\tau_{\hbar}{V}\neq0$ we have $\mathcal{H}^\mp_\Def=\ker(\d_V^\pm)$ by Mellin transform and the identity theorem. An interesting case is to consider potentials $(1-\tau_\Def)\Delta_4\Psi$ and $(1-\tau_\Def)H_4\Delta_4\Psi\in\mathcal{H}_\Def^-$, here we have a to \ref{QuoEigen2} quite analogous spectral realization.

The $\Conv$ derivation $\Mu_{\ln(t)}$ is not well-defined for $\mathcal{H}^\lambda_\pm(\Def)$ and $\mathcal{H}/\langle\widehat{\mathcal{Z}}^\lambda_\pm\rangle$, nevertheless $\Mu_{\ln(t)}$ induces an operation $H(\mathcal{H},\d^\pm_V)\rightarrow{H}(\mathcal{H},\d^\mp_V)$:\begin{lemma}\label{CohomologyDerivation}
Assume $\partial$ is a derivation of $\Conv$ and satisfies $\partial\tau=-\tau\partial$. If $\mathcal{R}^\mp=\ker(\d_V^\pm)$ we have a well-defined map $\partial:H(\mathcal{R},\d^\pm_V)\rightarrow{H}(\mathcal{R},\d^\mp_V)$ by $\partial\lfloor{f}\rfloor_V^\pm:= \lfloor{\partial{f}}\rfloor_V^\mp$.
\end{lemma}
 
\begin{proof}[Proof]  The assumptions imply the computation
$\d^\mp_V\partial{f}=\partial\d^\pm_V{f}-\left((\id+\tau)\partial{V})\right)\Conv\left((\id\pm\tau)f\right)$.
\end{proof}
Inspired by the strategy proposed in \cite{Den} it seems natural to consider for $\rho\geq0$ the \textit{\textbf{Brujin-Newman-P\'olya}} operators
$e^{-\rho\Mu_{\ln(t)}^2}=\sum_{n=0}^\infty\frac{-\rho^n}{n!}\Mu_{\ln^{2n}(t)}:H(\mathcal{H},\d^\pm_V)\rightarrow{H}(\mathcal{H},\d^\pm_V)$. 
\begin{lemma}\label{CohomologyPS}
The identities  $\mathcal{Z}^\lambda_{\pm'}(\Def)\d^\pm_V=\d^{\pm\pm'}_V\mathcal{Z}^\lambda_{\pm'}(\Def)$ hold, hence $\mathcal{Z}^\lambda_{\pm'}(\Def):H(\mathcal{H},\d^\pm_V)\rightarrow{H}(\mathcal{H},\d^{\pm\pm'}_V)$ are well-defined operators by $\mathcal{Z}^\lambda_\pm(\Def)\lfloor{f}\rfloor_V^\pm:=\lfloor\mathcal{Z}^\lambda_\pm{f}\rfloor_V^\pm$. We have the compatibilities $\mathcal{Z}^\lambda_+[\lfloor{f_1}\rfloor_V^+,\lfloor{f_2}\rfloor_V^+]^\ast=[\mathcal{Z}^\lambda_+\lfloor{f_1}\rfloor_V^+,\lfloor{f_2}\rfloor_V^+]^\ast$ and $\mathcal{Z}^\lambda_+\lfloor{f_1}\rfloor_V^-\Conv\lfloor{f_2}\rfloor_V^-=(\mathcal{Z}^\lambda_+\lfloor{f_1}\rfloor_V^-)\Conv\lfloor{f_2}\rfloor_V^-$.
\end{lemma}

\vfill{\begin{center}
\textcircled{c} Copyright by Johannes L\"offler, 2014, All Rights Reserved
\end{center}}


\subsubsection*{\textit{References}}
\begin{thebibliography}
 \renewcommand{\labelenumi}{[\theenumi]}
\begin{enumerate}
\vspace{-0.2cm}\bibitem{BK}
M. Berry, J. Keating,
\newblock $H=xp$ and the Riemann zeros,
\newblock Supersymmetry and Trace Formulae: Chaos and Disorder, New York: Plenum, 355-367, (1999)
\vspace{-0.2cm}\bibitem{Bru}
N.G. de Brujin,
\newblock The Roots of Trigonometric Integrals
\newblock Duke Math. J. \textbf{17} (1), 197-226, (1950)
\vspace{-0.2cm}\bibitem{Bur}
J.-F. Burnol,
\newblock The explicit formula in simple terms,
\newblock eprint:arxiv math.NT/9810169,v2 \textbf{22} Nov.(1998)
\vspace{-0.2cm}\bibitem{Co}
A. Connes,
\newblock  Noncommutative geometry and the Riemann zeta function,
\newblock Mathematics: Frontiers and Perspectives \textbf{2000}, V. Arnold et.al., eds., 35-55, (2000)
\vspace{-0.2cm}\bibitem{Co1}
A. Connes
\newblock Formule de trace en G\'eom\'etrie non commutative et hypoth\`ese
de Riemann,
\newblock C. R. Acad. Sci. Paris Ser. I Math. \textbf{323} n° 12, 1231-1236.(1996)
\vspace{-0.2cm}\bibitem{Co2}
 A. Connes, 
\newblock  Trace formula in noncommutative,
geometry and the zeros of the Riemann zeta function,
\newblock Selecta Math. (N.S.) \textbf{5} no. 1, 29-106,(1999)
\vspace{-0.2cm}\bibitem{ConMar}
A. Connes, M. Marcolli
\newblock  Noncommutative geometry, quantum fields and motives,
\newblock American Mathematical Society, Colloquium Publications, (2008)
\vspace{-0.2cm}\bibitem{Cs}
G. Csordas, T. S. Norfolk, R. S. Varga,
\newblock  The Riemann hypothesis and the Tur\'an inequalities,
\newblock Transactions of the American Mathematical Society, Volume \textbf{296}. Number 2, 521-541, (1986)
\vspace{-0.2cm}\bibitem{Den}
C. Deninger,  W. Singhof,
\newblock Real polarizable Hodge structures arising from foliations,
\newblock Annals of Global Analysis and Geometry \textbf{21}, 377-399, (2002)
\vspace{-0.2cm}\bibitem{KIMLEE}
H. Ki, Y-O. Kim, J. Lee,
\newblock  On the de Brujin-Newman constant,
\newblock Adv. Math. \textbf{222}, no. 1, 281-306, (2009)
\vspace{-0.2cm}\bibitem{Lag}
J. Lagarias,
\newblock  Zero spacing distributions for differenced $L$-functions,
\newblock Acta Arithmetica \textbf{120}, No. 2, 159-184, (2005)
\vspace{-0.2cm}\bibitem{RM}
R. Meyer
\newblock A spectral interpretation for the zeros of the Riemann zeta function,
\newblock Mathematisches Institut, Georg-August-Universit\"at G\"ottingen : Seminars Winter Term 2004/2005, Universit\"atsdrucke G\"ottingen, 117-137, (2005)
\vspace{-0.2cm}\bibitem{New}
C.M. Newman,
\newblock  Fourier transforms with only real zeros,
\newblock Proc. Amer. Math. Soc. \textbf{61}, 245-251, (1976)
\vspace{-0.2cm}\bibitem{Odl}
A.M. Odlyzko,
\newblock  An improved bound for the de Brujin-Newman constant,
\newblock Numerical Algorithms \textbf{25}, 293-303, (2000)
\vspace{-0.2cm}\bibitem{PO}
G. P\'olya,
\newblock  \"Uber die funktionentheoretischen Untersuchungen von J.L.W.V. Jensen,
\newblock Kgl. Danske Vid. Sel. Math-Fys. Medd. \textbf{7}, 3-33, (1927)
\vspace{-0.2cm}\bibitem{PO1}
G. P\'olya,
\newblock  \"Uber trigonometrische Integrale mit nur reelen Nullstellen,
\newblock Journal f\"ur die reine angewandte Mathematik, vol. \textbf{158}, 6-18, (1927)
\vspace{-0.2cm}\bibitem{RI}
B. Riemann,
\newblock \"Uber die Anzahl der Primzahlen unter einer gegebenen Gr\"osse,
\newblock Gesammelte Werke, Teubner (1892)
\vspace{-0.2cm}\bibitem{Sp}
A. Speiser,
\newblock Geometrisches zur Riemannschen Zeta Funktion,
\newblock Math. Annale. \textbf{110}, 514-521, (1934)
\vspace{-0.2cm}
\bibitem{Su1}
M. Suzuki,
\newblock  Nearest neighbor spacing distributions for the zeros of the real or imaginary part of the Riemann xi-function on vertical lines,
\newblock Acta Arith. 170,  no.  \textbf{1}, 47-65, (2015)
\vspace{-0.2cm}\bibitem{Ti}
E. C. Titchmarsh,
\newblock The theory of the Riemann zeta function,
\newblock The Clarendon press Oxford university (1986)
\vspace{-0.2cm}\bibitem{We}
A. Weil,
\newblock Fonctions zeta et distributions,
\newblock Seminar Bourbaki No. \textbf{312}, 158-163, (1966)
\end{enumerate}
\end{thebibliography}
\end{document}